\documentclass[12pt]{article}
\usepackage[left=2.3cm,right=2.3cm,top=2cm,bottom=2.5cm]{geometry}
\usepackage{amsfonts}
\usepackage{amsmath}
\usepackage{amsthm}
\usepackage{hyperref}
\usepackage{breqn}
\usepackage[usenames, dvipsnames]{color}
\allowdisplaybreaks 
\usepackage{graphicx}
\usepackage{caption}
\usepackage{subcaption}
\usepackage{float}
\usepackage{mathtools}
\usepackage{accents}
\usepackage{amssymb}
\usepackage{dsfont} % for unit operator  (indicator function) \mathds{1}

\providecommand{\keywords}[1]
{
	\small	
	\textbf{Keywords:} #1
}

\newcommand{\be}{\begin{equation}}
\newcommand{\bel}[1]{\begin{equation}\label{#1}}
\newcommand{\ee}{\end{equation}}
\newcommand{\bea}{\begin{eqnarray}}
\newcommand{\eea}{\end{eqnarray}}
\newcommand{\beann}{\begin{eqnarray*}}
\newcommand{\eeann}{\end{eqnarray*}}
\newcommand{\balign}{\begin{align}}
\newcommand{\ealign}{\end{align}}
\newcommand{\ba}{\begin{array}}
\newcommand{\ea}{\end{array}}

\newcommand{\bra}[1]{\mbox{$\langle \, {#1}\, |$}}
\newcommand{\ket}[1]{\mbox{$| \, {#1}\, \rangle$}}
\newcommand{\exval}[1]{\mbox{$\langle \, {#1}\, \rangle$}}
\newcommand{\inprod}[2]{\mbox{$\langle \, {#1} \, | \, {#2} \, \rangle$}}

\newcommand{\bfeta}{\boldsymbol{\eta}}

\newcommand{\R}{{\mathbb R}}
\renewcommand{\S}{{\mathbb S}}

\newcommand{\Z}{{\mathbb Z}}
\newcommand{\N}{{\mathbb N}}
\newcommand{\T}{{\mathbb T}} % e.g. for torus

\newtheorem{theorem}{Theorem}[section]
\newtheorem*{theorem*}{Theorem}
\newtheorem{definition}{Definition}[section]
\newtheorem{corollary}{Corollary}[section]
\newtheorem{lemma}{Lemma}[section]
\newtheorem*{lemma*}{Lemma}
\newtheorem{remark}{Remark}[section]
\newtheorem{proposition}{Proposition}[section]

\usepackage[dvipsnames]{xcolor}

\allowdisplaybreaks
\begin{document}

\newcommand{\nop}{\circ} \newcommand{\yesp}{\bullet}

%opening
\title{Open interacting particle systems and Ising measures}

\author{Ngo P.N. Ngoc\footnote{Institute of Research and Development, Duy Tan University, Da Nang 550000, Vietnam and Faculty of Natural Sciences, Duy Tan University, Da Nang 550000, Vietnam. email: ngopnguyenngoc@duytan.edu.vn} \and Gunter M. Sch\"{u}tz\footnote{
Departamento de Matem\'atica, Instituto Superior T\'ecnico, Universidade de Lisboa, Av. Rovisco Pais 1, 1049-001 Lisbon, Portugal. email: gunter.schuetz@tecnico.ulisboa.pt
%IAS 2, Forschungszentrum J\"ulich, 52425 J\"ulich, Germany. email: g.schuetz@fz-juelich.de
} 
}

\date{30/01/2025}

\maketitle

\begin{abstract}
We first survey some open questions concerning stochastic interacting particle systems with open boundaries. Then an asymmetric exclusion process with open boundaries that generalizes the lattice gas model of Katz, Lebowitz and Spohn (KLS) is introduced and invariance of the one-dimensional Ising measure is proved. The stationary current is computed in explicit form and is shown to exhibit
current reversal at some density. Based on the extremal-current principle for one-dimensional driven diffusive systems with one conservation law, the phase diagram for boundary-induced phase transitions is conjectured for this case. There are two extremal-current phases, unlike in the 
conventional open
asymmetric simple exclusion process which exhibits only one extremal-current phase or the previously 
considered conventional open KLS model with one or three extremal-current phases.
\end{abstract}


\section{Introduction}

Ever since the seminal paper by Spitzer on interacting random walks in 1970 \cite{Spit70} the study of stochastic interacting particle systems (SIPS), i.e., Markov processes describing the
random motion of particles on a lattice \cite{Ligg99}, has attracted a tremendous amount of attention in statistical physics
and in mathematical probability theory. Both rigorous 
and non-rigorous analytical and numerical work on such 
in general nonreversible systems has
shaped our understanding of nonequilibrium steady states, the emergence of hydrodynamic equations, fluctuations, and large deviations in 
simplified models of liquids and other many-body systems out of thermal equilibrium.

In particular, it has turned out that unlike in thermal equilibrium the boundaries of such 
systems,  usually defined on a subset of $\Z^d$, may determine in some substantial manner the bulk behaviour even when interactions are short-ranged. The 
fundamental reason for this phenomenon is the transport of conserved quantities such 
as mass or energy that determines how the bulk of a big system responds to boundary conditions that regulate the exchange of these quantities 
with the environment in which the system under consideration is embedded and hence controls the inflow and outflow of conserved quantities
from and to external reservoirs by a macroscopic amount. 

Transport in open systems may arise purely as a consequence of  gradients between different reservoir densities of the conserved quantities, in which case we speak of boundary-driven systems, or 
in the presence of further external bulk driving forces. Then we use the
term open bulk-driven systems.
Either way, 
the interdependence of bulk and boundary
has particularly dramatic consequences 
in one dimension. Transport may be anomalous and violate Fourier's or Fick's law \cite{Lepr16} and stationary distributions may exhibit long-range correlations \cite{Spoh83} or boundary-induced phase transitions \cite{Krug91}. It is therefore no surprise that particular attention has been paid to the 
study of one-dimensional open SIPS that are in contact with reservoirs at their two
boundaries. 
%Additional motivation to consider one-dimensional systems
%comes from numerous applications not only in physical systems
%but also in biological processes and other complex systems \cite{Chow24}.

The paradigmatic example of a SIPS is the asymmetric simple exclusion process
(ASEP) in which particles jump randomly on an integer lattice from one lattice site to a neighboring site with some spatial bias
and subject to the exclusion principle which allows for at most one particle on each site. The ASEP, defined precisely below,  was first invented as a mathematical model for protein synthesis by MacDonald, Gibbs, and Pipkin 
in 1968 \cite{MacD68} and was subsequently (and apparently independently) introduced by Spitzer in the
mathematical literature in his 1970 paper \cite{Spit70}. Since then applications of this model
and many generalizations of it have gone much beyond the realms of physics and mathematics
and include complex systems ranging from molecular motors in biological
cells to vehicular traffic on automobile highways \cite{Chow24}. Indeed, often the ASEP
is dubbed the Ising model of nonequilibrium physics since it is just as simple but mathematically and physically profound
as the celebrated Ising model is for equilibrium systems.
It has become a playground for the exploration of fundamental properties of quite general one-dimensional particle systems.

Mathematical challenges in the consideration of spatially finite
open systems (as opposed to their periodic or infinite counterparts) arise 
from a breakdown of
translation invariance, the loss of strict conservation of the bulk-conserved
quantities, and possibly the loss of other symmetries due to the presence
of boundary terms in the generator of the SIPS. Often specific tools must be 
devised to extract information on specific models.
It is not the purpose of this short article to survey the vast repertoire of methods
to treat open systems and even less to review the immense body of results
obtained by these methods. 
Instead, the goal is rather modest.

First, the focus is completely on steady states, i.e., invariant measures of SIPS. In Sec. 2 an attempt is made to give some structure to the
nature of the questions to be asked and to outline some open problems 
regarding the emergence of one-dimensional non-equilibrium phenomena.
In Sec. 3 we define the generator of a fairly general class of exclusion
processes for which these questions can be addressed. Deterministic large-scale behaviour under
Eulerian scaling, fluctuations on diffusive and superdiffusive level,
and large deviations are only
mentioned with a handful of references from which the reader can glean more information. 
The reader may also consult the older monographs and articles \cite{Ligg85,Derr98,Kipn99,Ligg99,Schu01,Bert15}
in which the ASEP features prominently. 

Second, marking an important anniversary in 2025, viz. the solution by Ernst Ising 
in 1925 of the one-dimensional Ising model \cite{Isin25},
a small piece of progress in the study of generalized asymmetric
simple exclusion processes is reported in Sec. 4. Originally, the Ising model was conceived by Ising's academic teacher Wilhelm Lenz as an equilibrium model for ferromagnetism. Here we prove invariance of the equilibrium measure
of the one-dimensional Ising model w.r.t. the Markovian dynamics of a generalization of the nonequilibrium
lattice gas model of Katz, Lebowitz, and Spohn for fast ionic conductors \cite{Katz84}. This result allows for computing the
exact stationary current from which the phase diagram of boundary-induced 
phase transitions is derived non-rigorously by appealing to the extremal-current
principle
for one-dimensional SIPS with one conservation law \cite{Popk99}. 

Our approach to defining boundary processes is constructive
and allows for constructing quite general lattice gas models with open boundaries and an 
explicit invariant measure by generalizing the family of generators 
introduced in Sec. 3. Thus it can be applied e.g.
for modelling initiation and termination of RNA or protein 
synthesis by molecular motors such as ribosomes or RNA
polymerase \cite{Chou11,Chow24} even when not only ASEP dynamics for translocation
of the motors is considered but also the internal mechanochemical cycle of the motor \cite{Trip08,Ngoc25}.
Also an extension to SIPS with long-range interactions
or more than one conserved quantity
looks feasible.

The paper ends with an appendix in which some important properties of the
one-dimensional Ising measure and their derivation using the transfer matrix method are summarized. This is textbook material that can be found in many places, but
for self-containedness, to adapt to our use of notation, and since the transfer matrix method appears to be not very well-known in the probabilistic literature, some readers might find this appendix useful.

\section{Open problems with open boundaries}

Despite the original motivation for the introduction of SIPS, viz., the study of large scale behaviour of microscopic models, many remarkable phenomena have been
revealed through the study of nonequilibrium steady states on the  microscopic level itself,
in particular, by numerical simulation and by explicit exact computation of invariant measures.
Yet, after several decades of work on the topic there are still many unresolved questions, both of general nature and with regard to specific models. 
Each of the very brief discussions of some of these questions
presented below
serves to outline open questions for future work.

\subsection{The setting}
\label{Sec:setting}

To fix some of our terminology which is motivated by a statistical physics
perspective and to illustrate the ideas that are presented more generally below, we first describe informally the ASEP
with open boundaries \cite{Ligg75,Krug91,Derr93a,Schu93b,Sand94b}.
A precise definition in terms of the generator is given in Sec. 3.

\begin{itemize}
\item Each lattice site of the finite integer lattice $\{1,\dots,L\}$ is either empty
or occupied by at most one particle.
\item Particles on the bulk sites $k\in\{2,\dots,L-1\}$ attempt to jump 
to right (left) neighbouring site $k+1$ ($k-1$). If
the target site $k\pm 1$ is empty, the jump succeeds, 
otherwise the jump attempt is rejected and the particle remains on site $k$.
Jump attempts occur independently of each other after exponentially
distributed random times with configuration-dependent rates
that we denote by $r_k(\bfeta)= r\eta_k(1-\eta_{k+1})$ for jumps from site $k$ to the right and 
$\ell_{k+1}(\bfeta)=\ell (1-\eta_{k})\eta_{k+1}$ for jumps from site $k+1$ to the left.
\item A particle on the left boundary site $1$ jumps with rate $r_1(\bfeta)$ to the
right (if site $2$ is empty) or is removed with rate $\gamma(\bfeta) = \gamma (1-\eta_1)$.
If site $1$ is empty then a particle is inserted with rate $\alpha(\bfeta) = \alpha \eta_1$.
\item A particle jumps on the right boundary site $L$ jumps with rate $\ell_L(\bfeta)$ to the
left (if site $L-1$ is empty) or is removed with rate $\beta(\bfeta) = \beta (1-\eta_L)$.
If site $L$ is empty then a particle is inserted with rate $\delta(\bfeta) = \delta \eta_L$.
\end{itemize}
In the absence of further specification of the jump rates we refer to this jump dynamics as simple exclusion
process (SEP) with open boundaries or closed boundaries if
$\alpha=\beta=\gamma=\delta=0$, i.e., when there is no exchange
of particles with the reservoir. 
When $r=\ell$ we speak of the symmetric simple exclusion process (SSEP)
and for $r\neq\ell$ the process is called asymmetric simple exclusion process (ASEP).
Specifically, when $\ell=0$ or $r=0$ the term totally asymmetric
simple exclusion process (TASEP) is used. 

More generally,
throughout the discussion of problems arising from open boundaries
we have in mind systems of particles with or without exclusion that are located on a finite
one-dimensional lattice with $d\geq 1$ lanes (such as $\{1,\dots,L\}^d$).
The term ``one-dimensional'' then means that we keep $d$ fixed, but allow for
the system length $L$ to become arbitrarily large. Second, we assume
short-range interactions between particles
by which we mean that transitions between microscopic particle configurations
involve only particles in a lattice segment with a fixed maximal number of
sites and rates for such a transition that either depend only on the particle configuration
in the vicinity of this segment or whose dependence on the location of particles far from this segment decays (in some sense that we do not qualify) sufficiently
quickly with the distance. The conserved quantities which may have an interpretation
mass, energy, or other physical quantities, are all represented by ``particles'' of different kinds, i.e., by discrete countable states in 
which each 
lattice site may be found. Conservation thus means that the transitions described above do not change 
the total number of each kind of particle, but only amount to a reshuffling
in space of some amount of the conserved quantities.
Each conservation law is represented by
exactly one type of particle, i.e., $n$ conserved quantities are represented
by $n$ distiguishable species of particles. Particles of a given species are,
however, indistiguishable.

Physically, the exclusion rule describes
an excluded-volume interaction (a hard-core repulsion) between identical particles and the bias
may be thought of as resulting either from some external driving force acting on the particles
or from some internal complex dynamics of e.g. microswimmers \cite{Elge15}, biological molecular motors \cite{Chou11,Chow24}, or other types of active particles \cite{Bowi22}. For closed boundaries only particle jumps occur and the process
is conservative, i.e., the total number of particles is conserved.
Consequently, on a finite lattice any such SIPS is trivially nonergodic in the sense that it has a family of invariant measures that is specified by the total number of particles on the lattice.

The sites near the left (right) edge of the lattice, i.e.,
near site 1 ($L$), are considered to be left (right) boundary sites. On these sites
transitions may take place that violate the conservation law, i.e., particles may be
inserted or removed and usually the invariant measure is unique. This describes the particle exchange of particles on the lattice with two virtual external reservoirs, left and right, whose states are not considered  and are physically only assumed to arise from some fixed reservoir
densities of the various particle species.
When there is no insertion of particles we speak of absorbing
boundaries. The invariant measure is then the trivial Dirac measure concentrated
on the configuration representing the empty lattice. 

Thus within this description, the SIPSs under consideration have a countable state space even in the absence of exclusion. Each state is called a configuration of particles. They live on a finite lattice, have short-range interactions, and describe the Markovian stochastic dynamics of $n$ species of particles whose total random number is 
countable at any given instant of time. To exclude non-generic problems we assume
the Markov chain to be ergodic and hence have a unique invariant measure.

Due to the bulk conservation law there is a particle current associated with each of the $n$ particle species. If the invariant measure is characterized by non-zero currents
then it represents a nonequilibrium
steady state (NESS) \cite{Zia06,deCa17}. Like the $n$ stationary particle densities these 
$n$ stationary
currents are (in general) complicated functions of system size and the transition 
rates.
In the thermodynamic limit $L\to\infty$ these $n$ currents
can be expressed as functions of the stationary particle densities.
Nonzero limiting currents arise from bulk driving forces that act on the particles
which are represented in the transition rates in some bias to the right (in positive
direction) or to the left (in negative direction). Conversely, the absence of currents in the thermodynamic limit $L\to\infty$ indicates that
locally the system is in state of thermal equilibrium. In a finite system
a current can then be maintained by different reservoir densities of the
particle species, i.e., in boundary-driven systems.

\subsection{Boundary-driven systems}

We discuss some issues related to the stationary behaviour
in the presence of boundary gradients in systems that
would be in equilibrium, i.e., reversible, in the absence of such gradients.

\paragraph{$\bullet$  Long-range correlations:}
Translation-invariant one-dimensional SIPS are generally not expected to exhibit long-range correlations \cite{Garr90}. Indeed, looking only
at the approach to stationarity of the local particle density of the boundary-driven SSEP which
is governed simply by the diffusion equation with Dirichlet, Neumann, or Robin boundary conditions
\cite{Bald17} 
as if the particles were non-interacting, one may be tempted
to conjecture that stationary long-range
correlations would indeed not exist.
 It is therefore exciting to note
that long-range correlations are actually expected 
to be a {\it general} feature of one-dimensional SIPS
when they are boundary driven. Such correlations
were first obtained for the open SSEP by duality \cite{Spoh83} at the level of the centered 
two-point density correlation function
and later studied for higher-order correlations in great detail
using the matrix product ansatz \cite{Blyt07}, originally developed in 
\cite{Derr93a} for the computation of the invariant measure of the open ASEP. 
The cumulant generating function 
yields an amplitude of order $L^{-n+1}$  for $n$-point truncated correlations \cite{Derr07} and the $n$-point correlations of the centered occupation numbers
appear to be given for even $n$ by the Wick theorem in terms of two-point correlations (which are of order $L^{-1}$) and by a generalized Wick theorem
in terms of two-point correlations and a three-point correlation (which is of order $L^{-2}$) \cite{Gonc25}. 

Quite generally, these correlations as well as the stationary distribution of the particle current can be described in the framework of large-deviation theory. The large-deviation
properties have been derived both using the matrix product ansatz in combination with an additivity principle \cite{Derr07a,Derr21,Saha25} and macroscopic fluctuation theory
\cite{Bert15,FranT23}. A problem of great interest would be to investigate
whether (or to which extent) these higher-order correlation properties are special to the SSEP or a universal feature of
boundary-driven systems. This is of interest also 
from an experimental perspective as the SSEP and its two-species
generalization have served successfully as models for
for interface growth in the Edwards-Wilkinson universality class
\cite{Edwa82,Krug97}, for  single-diffusion  \cite{Kukl96,Wei00} and for reptation dynamics %\cite{Doi85} 
of entangled
DNA \cite{Perk94,Schu99}. 

\paragraph{$\bullet$ Duality:}
Duality \cite{Ligg85,Jans14} is a powerful tool in the study of SIPS
as it relates time-dependent expectations of one process in terms
of possibly simpler time-dependent expectations of the dual process.
The dual open SSEP has absorbing boundaries \cite{Schw77,Cari13}, a property
that appears also in other boundary-driven systems 
\cite{Fras22,FranC23,Casi24} but does not seem to have a very well understood level of generality.
Specifically, it has not only revealed the long-range correlations in the SSEP
mentioned above but also an intimate connection of the stationary probabilities of the SSEP to absorption probabilities in the dual process
\cite{Flor23}. 
A class of models where further mathematically rigorous progress in this direction is a promising prospect includes symmetric exclusion of rods, i.e., of particles that cover more than one lattice site but jump by only one site.
Such SIPS have been introduced for
the study of biological molecular motors as these are usually larger than the sites on the macromolecular template along which they move \cite{MacD68,Laka03,Shaw03}. This model, for rods of size two equivalent to the ergodic sector of the facilitated exclusion process \cite{Barr23,Erig23}, can be mapped to the SSEP \cite{Scho04} and thus exhibits the same $SU(2)$ symmetry that gives rise to duality in the SSEP \cite{Schu94,Giar09}.

\paragraph{$\bullet$ Multi-species systems:} Single-diffusion occurs
when particles of different species move on a one-dimensional substrate
without being able to pass each other. In this case a single particle
moving in an environment of other particles moves subdiffusively
\cite{Arra83}, giving rise to non-standard hydrodynamic behaviour
on macroscopic scale \cite{Quas92}. For open boundaries this leads
to boundary-induced bulk phase transition
and a violation of Fick's law by allowing for partial uphill diffusion, i.e.,
a particle current against the boundary gradient \cite{Brza06a,Brza06b}.
Remarkably, in a specific two-species system without single-file dynamics it was noted that partial uphill diffusion is possible for the discrete particle system on the lattice, whereas it is lost in the hydrodynamic limit \cite{Casi23}. A deeper understanding of the role
of the single-file constraint for phase transitions and uphill diffusion, however, is still lacking.

\paragraph{$\bullet$ Long-range interactions:}

Also long-range interactions
play a role in applications. They arise either when the
particle jumps are long-ranged \cite{Frei02,Jara08,Bald17,Bern23},
in which case a fractional Laplacian and fractional
boundary conditions
may arise to describe both the steady state and its approach to it \cite{Gonc23a}, or when the particle jumps are short-ranged, but the rate of jumps does not decay sufficiently quickly with the distance of the other particles.
The latter case is capable of describing a large variety
of physical scenarios such as the large scale dynamics of the particle density following the
porous medium equation \cite{Card23,Gonc23b},
logarithmic interactions that appear in the dynamics of
DNA melting \cite{Bar09,Hirs11} and which may lead to
phase separation \cite{Beli24}, or exhibit a conformally
invariant hyperuniform phase \cite{Kare17}.

A general open question for such systems is a natural definition of reservoir dynamics with boundary gradients. 
The approach taken in this paper in Sec. \ref{Sec:KLS_model} for short-range interactions without boundary gradients
is constructive and may serve as a guide.
One may then be able to address the question of long-range
correlations due to boundary gradients.

\paragraph{$\bullet$  Open quantum systems:} Remarkably, for closed boundaries or periodic boundary
conditions the generator of the SSEP is identical
to the quantum Hamiltonian of the spin-1/2 Heisenberg ferromagnet \cite{Alex78}
and features reminiscent of the open SSEP were discovered also in the boundary-driven quantum case \cite{Pros11,Kare13b,Pros15,Buca16}. It would be important to understand whether the long-range correlations both in classical Markovian and quantum systems have universal properties.

\subsection{Bulk-driven open systems}

Many of the approaches used for boundary-driven systems
have been employed also for open bulk-driven systems.
In particular, the MPA can be applied without modification
and macroscopic fluctuation theory can be extended to account for weakly bulk-driven dynamics. A basic insight that has emerged in the study
of the effect of open boundaries already in the 1990s that was mentioned already above, viz., the appearance of boundary-induced phase transitions \cite{Krug91}, result from interplay of shocks and rarefaction waves when they reach the boundary \cite{Popk99}. They result in two distinct types of non-equilibrium
phase transitions: (i) First order, driven by domain-wall dynamics \cite{Kolo98b} where the fluctuations of the current become singular
\cite{Derr94}
and (ii) second-order to extremal-current phases which exhibit
a universal density profile \cite{Krug91,Oerd98,Hage01} and algebraically decaying correlations \cite{Derr93c}.
Given the stationary current for the infinite system, the phase diagram
can be obtained through an extremal current principle \cite{Popk99}
\be
\label{ex}
 j = \left\{ \ba{l} \displaystyle
\max_{\rho \in [\rho_+,\rho_-]} j(\rho) \quad \mbox{ for } \rho_->\rho_+ \\
\displaystyle
\min_{\rho \in [\rho_-,\rho_+]} j(\rho) \quad \mbox{ for } \rho_-<\rho_+.
\ea \right.
\ee
which predicts the selected stationary density of an open bulk-driven system
with different reservoir densities $\rho_{-}$ at the left boundary
and $\rho_{+}$ at the right boundary.

Rarefaction waves arise macroscopically as entropy solutions 
of the macrosocpic hyperbolic conservation law \cite{Gart88,Reza91}
obtained as hydrodynamic limit under Eulerian scaling \cite{Kipn99}.
Shocks, which in the original biological context of the ASEP correspond to ``traffic jams'' of ribosomes
\cite{MacD68,Schu97a,Chou11,Chow24}, are readily understood in terms of the macrosopic Rankine-Huginot condition for the stability of shocks \cite{Lax73}. 
Both phenomena can be derived microscopically by studying the flow
of small perturbartions \cite{Popk99}.
A phenemenological discussion appealing to non-linear
fluctuation hydrodynamics \cite{Krug91}, renormalization group
and scaling arguments \cite{Oerd98,Hage01}, and also mathematically rigorous results
for attractive particle systems \cite{Baha12} lend support to the picture
developed in \cite{Popk99} for open boundaries which is discussed in some detail in 
\cite{Schu01}. 
Some fundamental problems, however, have remained unresolved.

\paragraph{$\bullet$ Boundary-induced phase transitions:}
A fundamental argument used in \cite{Popk99} to explain 
boundary-induced
phase transition of first order is the notion that 
macroscopic shock discontinuities that satisfy the Rankine-Huginot condition
are also microscopically sharp, as proved e.g. for the ASEP in
\cite{Ferr91,Derr93b,Beli02}. It is not clear, however, what the
microscopic structure of a contact discontinuity looks like where the
Rankine-Huginot criterion does not hold. Numerical evidence
coming from the KLS-model with open boundaries
where such a contact discontinuity may arise \cite{Schw24} suggests that in a finite open system there still exists a sharp change of density that on Eulerian
scale would be a shock, but the intrinsic width of this macroscopic
discontinuity grows algebraically with system size, i.e., such a shock
is {\it not} microscopically sharp. In \cite{Schw24} a heuristic attempt is made
to relate this scaling to current fluctuations. However, a fully satisfactory
understanding of this scaling and, more profundly, of the microscopic
structure of contact discontinuities is lacking.

\paragraph{$\bullet$ Duality:}
Given the wealth of results on duality for boundary-driven systems,
it is somewhat surprising that insightful dualities for open bulk driven systems have remained rather elusive \cite{Ohku17}. An intriguing result
relates by duality a bulk-driven to a boundary driven process \cite{Cari24}
and it would be interesting to explore the level of generality of such
a construction. Duality for the open ASEP has been studied in \cite{Barr24}
in terms of a system of ordinary differential equations. This system is
solved by Bethe ansatz for the ASEP on the open integer lattice, i.e., with only one
reservoir. The extension to the finite lattice with open boundaries
remains an open problem. This problem was attacked also very recently
by means of reverse duality \cite{Schu23a,Schu23b}, which unlike conventional duality
relates the time-dependent probability measure of two Markov chains
rather than time-dependent expectations. Thus for a specific parameter
manifold of particle insertion and removal rates in the ASEP 
the picture of a microscopically
sharp shock that is crucial for understanding the first-order boundary-induced phase transition could be proved rigorously
and refined in terms of an exact random walk property
of shocks. The extension of reverse duality to other models and its relation
to the bulk symmetries of the generator of the process, however, needs to be clarified.

\paragraph{$\bullet$ Multi-species systems:}
In the presence of multiple bulk conservation laws new phenomena
become possible, including phase 
separation \cite{Evan98b,Maha20} and an intriguing dynamical universality classes of fluctuations \cite{Spoh14,Popk15b,Ahme22,Cann24,Popk24}, that have no analogue in thermal equilibrium or in the case of only one conserved quantity.
In open systems shocks still play a decisive role for understanding
boundary-induced phase transitions, but the picture becomes much more
complex and is only partially understood. The main novelty is the
reflection of shocks at the boundaries \cite{Popk04}
which results in both boundaries ``communicating'' with each other.
Hence they cannot be treated independently by means of a simple
extremal-current principle as in the case of one conservation law.
This results in rather complex phase diagrams of boundary-induced
phase transitions which is poorly understood. Significant progress
has been achieved recently for an integrable two-species exclusion
process \cite{Cant08,Cant22} which suggests a generalization of the
extremal current principle \eqref{ex} in terms of Riemann invariants
\cite{Popk04,Cant24}. Nevertheless, full understanding is still missing,
particularly in the case when bulk fluctuations belong to different
dynamical universality classes \cite{Rako04,Ahme22,Cann24}.

\paragraph{$\bullet$ Long-range interactions:}
Also for open bulk-driven systems with long-range interactions
a good general approach to defining reservoir dynamics with boundary 
gradients is an open question. This question arises e.g. for long-range
exclusion interaction with nearest-neighbor jumps where phase 
separation can occur in a conservative periodic system \cite{Priy16,Beli24}. In such a scenario it is unclear whether phase separation prevails for non-conservative open boundaries and, if it should be the case, whether the active phase is pinned to one of the boundaries
and whether the Gallavotti-Cohen symmetry of the particle current
breaks down as in the related phenomenon of condensation
in the open asymmetric zero-range process \cite{Harr06}.

\paragraph{$\bullet$ Open quantum systems:}
The shocks that feature so prominently in the ASEP with open boundaries
have a quantum analogue as spin helix states \cite{Popk17,Popk21}
which have recently been observed experimentally in Heisenberg quantum magnets \cite{Jeps22}. This is particularly remarkable as quantum
spin chains have been shown in recent years to display properties in the universality class of the Kardar-Parisi-Zhang equation \cite{Ljub19,Ye22,DeNa23} that 
characterizes also the classical ASEP \cite{Corw12,Spoh17,Mate21,Quas23,Arai24}.
However, little is known about the
dynamics of these helix states and whether this problem can be attacked by means of duality or reverse duality.

\section{Open simple exclusion processes}

Attacking the problem of gaining general insights by studying concrete models is certainly an ambitious task. Perhaps the most promising candidates for further study are exclusion processes that generalize the SEP described above, an important one being the Katz-Lebowitz-Spohn (KLS) model with next-nearest-neighbour interactions \cite{Katz84}.

It looks natural to abbreviate the term simple exclusion process by SEP with S for simple. However, 
as this acronym is also used for the symmetric simple exclusion processes (with S
for symmetric) here the shorthand EP for any exclusion process
is used. The standard symmetric simple exclusion process (with {\it only} exclusion and no other interaction) and nearest neighbor jumps 
is then denoted by SSEP and the corresponding (equally standard) version of this model
but with asymmetric, i.e., spatially biased, jump rates by ASEP.

To account for this setting with some degree of generality we define such an {\it open} EP for a one-dimensional lattice of $L$ sites, labelled
by integers from the set  $\Lambda_L:=\{1,\dots,L\}$, as follows. The site $k=1$ ($k=L$)
is called left (right) boundary site, all other sites are called bulk sites.
The {\it occupation 
numbers} $\eta_k\in\{0,1\}$ indicate whether site $k$ is occupied ($\eta_k=1$) or
empty ($\eta_k=0$). The symbol $\bfeta := (\eta_{1}, \eta_{2}, \dots, \eta_{L})$ 
denotes the configuration on the whole lattice. 
The state space of the process
is $\Omega_L := \{0,1\}^L$.

For a precise definition of the process in terms of the generator it is 
convenient to define the vacancy occupation variables 
\begin{equation}
\bar{\eta}_{k} := 1 - \eta_{k},
\label{baretadef}
\end{equation}
the {\it flipped} configuration $\bfeta^{k}$ defined by the occupation numbers
\bel{etaflipdef}
\eta^{l}_k =\eta_k + (1-2\eta_l) \delta_{l,k} = \left\{\ba{ll}
\bar{\eta}_{k} & \mbox{if } k=l \\
\eta_{k} & \mbox{else} ,
\ea \right.
\ee
the
{\it swapped} configuration $\bfeta^{l,l+1}$ defined by the occupation numbers
\bel{etaswapdef}
\eta^{l,l+1}_k = \eta_k + (\eta_{l+1}-\eta_{l})(\delta_{l,k}-\delta_{k,l+1})
= \left\{\ba{ll}
\eta_{l+1} & \mbox{if } k=l \\
\eta_{l} & \mbox{if } k = l+1 \\
\eta_{k} & \mbox{else}.
\ea \right.
\ee
%and the {\it shifted} configurations $\bfeta^{\pm n}$ defined by the shifted occupation numers
%\bel{etashift}
%\eta^{\pm n}_k = \eta_{k\pm n \mbox{\footnotesize \ mod} \ L}.
%\ee

Due to the exclusion interaction the configuration-dependent bond jump rates across a bond $k,k+1$ are taken to be of the form
\begin{equation}
w_k(\bfeta) := r_k(\bfeta) + \ell_{k+1}(\bfeta)
\label{bjr}
\end{equation}
with functions
$r_k:\Omega_L \to \R^{+}$ and  $\ell_{k+1}:\Omega_L \to \R^{+}$ generalizing those of the ASEP defined in Sec. \ref{Sec:setting}.
The generator for a particle jump across bond
$(k,k+1)$ is then given by
\begin{equation}
\mathcal{L}_k f(\bfeta) = w_k(\bfeta) [f(\bfeta^{k,k+1}) -  f(\bfeta)]
\label{jumpgenbulk}
\end{equation}
for measurable functions $f:\{0,1\}^L\to \R$. According to the informal
description the particle jumps are encoded in the bulk generator
\bel{genASEPbulk}
\mathcal{L}^{bulk} f(\bfeta) = \sum_{k=1}^{L-1} \mathcal{L}_k f(\bfeta).
\ee

The particle exchange with the reservoirs is described by the configuration-dependent boundary rates
\begin{equation}
b^{-}(\bfeta) := \alpha(\bfeta) + \gamma(\bfeta), \quad 
b^{+}(\bfeta) := \beta(\bfeta) + \delta(\bfeta).
\label{boundaryrates}
\end{equation}
Here $\alpha(\bfeta)$ ($\delta(\bfeta)$), taken to be proportional to $\bar{\eta}_{1}$ ($\bar{\eta}_{L}$), is the rate
of particle insertion at the boundary site 1 ($L$) and 
$\gamma(\bfeta)$ ($\beta(\bfeta)$), taken to be proportional to $\eta_{1}$ ($\eta_{L}$), is the rate
of particle removal at the boundary site 1 ($L$). This
leads to the
boundary generators
\begin{equation}
\mathcal{B}^{-} f(\bfeta) = b^{-}(\bfeta) [f(\bfeta^{1}) -  f(\bfeta)], \quad
\mathcal{B}^{+} f(\bfeta) = b^{+}(\bfeta) [f(\bfeta^{L}) -  f(\bfeta)].
\label{genboundary}
\end{equation}
The generator $\mathcal{L}$ of an open EP of the form described above 
is thus given by
\bel{genASEP}
\mathcal{L} = \mathcal{B}^{-} + \mathcal{L}^{bulk} + \mathcal{B}^{+}.
\ee
%When in addition to this condition particles jump from site $L$ to site $1$
%with rate $r_L(\bfeta)$ and from site $1$ to site $L$ with rate $\ell_1(\bfeta)$
%we speak of toroidal boundary conditions. The process is called translation
%invariant when $r_{k+1}(\bfeta) = r_k(\bfeta^{+1})$ and $\ell_{k-1}(\bfeta) = \ell_k(\bfeta^{-1})$ for all $k\in\Lambda_L$. 
Exclusion processes with a generator
of the form \eqref{genASEP} are generically ergodic
and not reversible. We denote expectations in the invariant measure
by $\exval{\cdot}_L$.
When $\alpha(\cdot) = \beta(\cdot) = \gamma(\cdot) = \delta(\cdot) = 0$
the total particle number
\begin{equation}
N(\bfeta) := \sum_{k=1}^{L} \eta_k 
\label{Ndef}
\end{equation}
 is conserved and we speak of closed boundaries.

For open boundaries
the stationary density $\rho(L)$ of the process is defined by the normalized expectation
\begin{equation}
\rho(L) := \frac{1}{L}  \exval{N}_L
\label{rhoLdef}
\end{equation}
of the total particle number. Due to particle exchange at the boundaries this quantity is a (usually complicated) function of all the parameters that specify the jump rates
and boundary insertion and removal rates. The thermodynamic limit is defined
by
\begin{equation}
\rho := \lim_{L\to\infty} \rho(L) .
\label{rhodef}
\end{equation}
The stationary density profile 
\begin{equation}
\rho_k(L) := \exval{\eta_k}_L
\label{rhokldef}
\end{equation} 
quantifies the
effect of the boundary reservoirs on the stationary local particle density.

Due to particle number conservation in the bulk, the action of the generator on the occupation number $\eta_k$ in the setting described above gives rise to the discrete microscopic
continuity equation
\begin{equation}
\mathcal{L} \eta_k = j_{k-1}(\bfeta) - j_k(\bfeta), \quad 1 < k < L
\end{equation}
with the  {\it instantaneous current}
\begin{equation}
j_k(\bfeta) := r_k(\bfeta) - \ell_{k+1}(\bfeta).
\label{jinstdef}
\end{equation}
Since by definition of invariance $\exval{\mathcal{L} \eta_k}_L = 0$ the stationary current 
\begin{equation}
j(L) := \exval{j_k}_L
\label{jLdef}
\end{equation}
is independent of $k$. The physical interpretation of this quantity,
denoted throughout this work by $j$,
is the stationary net flow of particles across a bulk lattice bond $(k,k+1)$.
It is a quantity of central importance in the study of nonequilibrium
steady states \cite{Derr98,Schu01} and hydrodynamic
limits \cite{Kipn99}. A nonzero stationary current indicates absence
of reversibility \cite{Zia06,deCa17}. 
Except when in the thermodynamic limit
\begin{equation}
j := \lim_{L\to\infty} j(L) 
\label{jdef}
\end{equation}
the limiting stationary current $j$ is zero for all densities $\rho$ the particle system
is called bulk-driven. Generically one expects for the bulk density profile
of  bulk-driven systems $\lim_{L\to\infty} \rho_{[xL]}(L)
= \rho$ for any $x\in(0,1)$ \cite{Krug91}. 

We stress that this set-up can be generalized in analogous fashion
for non-exclusion and multi-species processes, and to more complex
jump and boundary processes. In the following we consider a generalization to nearest-neighbor jumps
with next-nearest neighbor interactions and exclusion.

\section{Generalizations of the Katz-Lebowitz-Spohn 
model}\label{Sec:KLS_model}

From a physics perspective, the simple on-site exclusion interaction of the ASEP appears somewhat artificial and suited only to describe the interaction of hard spheres. More realistically, one has to consider short-range interactions. 

\subsection{The one-dimensional Katz-Lebowitz-Spohn 
model on the integer torus}
\label{KLS_model}

The simplest
generalization of this kind consists in postulating jump rates that depend not only on the 
target site of the jump, but also on the occupation numbers of the sites $k-1$, $k+2$
next to the bond $(k,k+1)$ across which a jump attempt occurs. The paradigmatic
model of this kind is the one-dimensional Katz-Lebowitz-Spohn (KLS) model \cite{Katz84} defined with the parameters $\kappa, \epsilon \in [-1,1]$
by the jump rates 
\bea
r_k(\bfeta) & = & r \eta_{k} \bar{\eta}_{k+1} \left[ 
(1+\kappa) (\bar{\eta}_{k-1} - \eta_{k+2}) + (1+\epsilon) \eta_{k-1}
+ (1-\epsilon) \eta_{k+2}\right]  
\label{rKLSbulkrate} \\
\ell_{k+1}(\bfeta) & = & \ell \bar{\eta}_{k} \eta_{k+1} \left[
(1+\kappa)(\bar{\eta}_{k-1} - \eta_{k+2}) +
(1 - \epsilon) \eta_{k-1} + (1+\epsilon) \eta_{k+2} 
\right] 
\label{lKLSbulkrate}
\eea
for all $k\in\{1,\dots,L\}$ where the site index $k$ of the
occupation variables is defined modulo $L$, i.e., $\eta_{k+m L} := \eta_k$ for
$k\in\Lambda_L$ and all $m\in\Z$.
The jump rates can be visualized by the pictorial presentation
\begin{align}
	& 0\accentset{\curvearrowright}{10}0 \xrightarrow{r(1+\kappa)} 0010,\quad 
	1\accentset{\curvearrowright}{10}0 \xrightarrow{r(1+\epsilon)} 1010,\quad 
	0\accentset{\curvearrowright}{10}1 \xrightarrow{r(1-\epsilon)} 0011, \quad  
	1\accentset{\curvearrowright}{10}1 \xrightarrow{r(1-\kappa)} 1011
\label{rKLS}\\[4mm]
	&0\accentset{\curvearrowleft}{01}0 \xrightarrow{\ell(1+\kappa)} 0100,\quad  
	1\accentset{\curvearrowleft}{01}0 \xrightarrow{\ell(1-\epsilon)} 1100,\quad  
	0\accentset{\curvearrowleft}{01}1 \xrightarrow{\ell(1+\epsilon)} 0101, \quad  
	1\accentset{\curvearrowleft}{01}1 \xrightarrow{\ell(1-\kappa)} 1101
\label{lKLS}
\end{align}
where the curved arrows indicate the jump direction. 

The generator for periodic boundary conditions on the integer torus
$\T_L$ is given by
\bel{perKLS}
\mathcal{L} f(\bfeta) = \sum_{k=1}^{L} \mathcal{L}_k f(\bfeta)
\ee
where $\mathcal{L}_k$ is of the form \eqref{jumpgenbulk} with jump rates $r_{k}$ and $\ell_{k+1}$ given by \eqref{rKLSbulkrate} and \eqref{lKLSbulkrate}.

As pointed out in \cite{Katz84} the one-dimensional Ising measure
\begin{equation}
	\pi_{L}^\ast(\boldsymbol{\eta}) = \dfrac{1}{Z_{L}} \pi_{L}(\boldsymbol{\eta})
\label{Isingperdef}
\end{equation}
with the Boltzmann weight
\begin{equation}
\pi_{L}(\boldsymbol{\eta}) =
\mathrm{e}^{-\beta \left[E(\boldsymbol{\eta}) - \varphi N(\boldsymbol{\eta})\right]} ,
\label{IsingBW}
\end{equation}
the static interaction energy 
\begin{equation}
E(\boldsymbol{\eta}) =
J \sum_{k=1}^{L}\eta_k\eta_{k+1} ,
\label{IsingE}
\end{equation}
and the partition function
\begin{equation}
Z_{L} = \sum_{\boldsymbol{\eta}\in\Omega_L} \pi_{L}(\boldsymbol{\eta})
\label{Zper}
\end{equation}
is invariant for all $\varphi\in\R$ if and only if 
\begin{equation}
\mathrm{e}^{-\beta J} = \frac{1-\epsilon}{1+\epsilon}.
\label{bJeps}
\end{equation}
The Ising measure is well-defined for $\varphi\to\pm \infty$
but reduces in these limits to the Dirac measure concentrated on the completely empty
or completely occupied configuration respectively. These trivial limiting cases are excluded from consideration. 
%Also the limits 
%$J\to\pm \infty$, corresponding to $\epsilon \to \pm 1$, are well-defined
%for finite system size $L$.
%For $\epsilon=1$,
%corresponding to the limit $J\to\infty$, the invariant measure \eqref{Isingperdef} is uniform on the set of configurations that
%have no neighboring particles, i.e., the set $\{\bfeta\in\Omega_L:
%\eta_{k} \eta_{k+1} = 0 \, \forall \, k \in \T_L\}$.
If $r=\ell$ the model is symmetric under space reflection $k\mapsto L+1-k$
and reversible. 

In equilibrium thermodynamics, the invariant measure \eqref{Isingperdef} is an Ising measure
for periodic boundary conditions where the non-negative real parameter $\beta$ is proportional to the inverse of the 
temperature $T$, i.e., $\beta = 1/(k_BT)$, where $k_B$ is the Boltzmann constant and $J\in\R$ is a static interaction strength. Since $\beta$ appears only on the product
$\beta J$ we set $\beta=1$.
The parameter $\varphi$ plays the role of chemical potential that regulates the particle density $\rho(L)$ through the relation
\begin{equation}
\rho(L) = - \frac{1}{L} \partial_\varphi \ln{Z_{L}}. 
\end{equation}
%As a function of the chemical potential $\varphi$ the density $\rho_{L}$ is strictly monotone and hence an invertible function of $\varphi$. 

Expectations w.r.t. the Ising measure \eqref{Isingperdef} can be computed in explicit form 
using standard transfer 
matrix techniques \cite{Baxt82}. As this approach appears to be not well-known
in a probabilistic setting the main steps are summarized in Appendix
\ref{App:Ising}. In particular, the particle density is given in explicit form by
\begin{equation}
\rho(L) = - \frac{1}{L} \partial_\varphi \ln{(\lambda_{0}^L + \lambda_{1}^L)}
\end{equation}
in terms of the eigenvalues $\lambda_{0,1}$ \eqref{eigenvalues} of the transfer matrix for the
Ising measure.
Expectations in the thermodynamic limit $L\to\infty$
are denoted by $\exval{\cdot}$ without the subscript $L$. These quantities deviate
from the corresponding
expectations for fixed $L$ by finite-size-corrections that are exponentially small 
in $L$ with parameter $\ln{(\lambda_{1}/\lambda_{0})}$. In particular, the density 
\begin{equation}
\rho  := \lim_{L\to\infty} \rho(L) = - \partial_\varphi  \ln{\lambda_{0}}
\end{equation}
is given by the largest eigenvalue $\lambda_{0}$ of the transfer matrix.

\subsection{Generalized periodic KLS model}

For any fixed $\epsilon \in [-1,1]$ the periodic Ising measure \eqref{Isingperdef} is invariant for all choices of the kinetic interaction parameter $\kappa \in [-1,1]$ and it is also independent of the hopping asymmetry
$q = \sqrt{r/\ell}$. Without changing the invariant measure, the KLS model can therefore be generalized to allow for different
kinetic interaction parameters for jumps to the right and to the left respectively. This simple fact leads to our first result which can be stated without further proof.

\begin{theorem}
\label{maintheory_1}
For all $\varphi\in\R$ the Ising measure \eqref{Isingperdef} with static interaction parameter 
given by \eqref{bJeps} is invariant for a generalized KLS model defined
on $\T_L$ with periodic boundaries 
by the generator 
\bel{genKLSper}
\mathcal{L} f(\bfeta) = \sum_{k=1}^{L}
w_k(\bfeta) [f(\bfeta^{kk+1}) -  f(\bfeta)]
\ee
 with jump rates $w_k(\bfeta) = r_k(\bfeta) + \ell_{k+1}(\bfeta)$ given by
\bea
r_k(\bfeta) 
& = & r \eta_{k} \bar{\eta}_{k+1} \left[ 
(1+\kappa) (\bar{\eta}_{k-1} - \eta_{k+2}) + (1+\epsilon) \eta_{k-1}
+ (1-\epsilon) \eta_{k+2}\right] 
\label{rGKLS} \\
& = & \eta_{k} \bar{\eta}_{k+1} \left[ 
r (1+\kappa)  + r (\epsilon-\kappa) \eta_{k-1}
- r (\epsilon+\kappa) \eta_{k+2}\right] 
\label{rGKLS2} \\
\ell_{k+1}(\bfeta) & = & \ell \bar{\eta}_{k} \eta_{k+1} \left[
(1+\lambda)(\bar{\eta}_{k-1} - \eta_{k+2}) +
(1 - \epsilon) \eta_{k-1} + (1+\epsilon) \eta_{k+2} 
\right]
\label{lGKLS}  \\
& = & \bar{\eta}_{k} \eta_{k+1} \left[
\ell (1+\lambda) - 
\ell (\epsilon+\lambda) \eta_{k-1} + \ell (\epsilon-\lambda) \eta_{k+2} 
\right]
\label{lGKLS2}  
\eea
for all $\epsilon,\kappa, \lambda \in [-1,1]$ and
for all $r,\ell \in \R_0^{+}$.
\end{theorem}

In pictorial presentation the generalized jump rates are like
those of the standard KLS model displayed in \eqref{rKLS},
\eqref{lKLS}, except that in \eqref{lKLS} the parameter $\kappa$
has to be replaced by $\lambda$.
%\begin{align*}
%	& 0\accentset{\curvearrowright}{10}0 \xrightarrow{r(1+\kappa)} 0010,\quad 
%	1\accentset{\curvearrowright}{10}0 \xrightarrow{r(1+\epsilon)} 1010,\quad 
%	0\accentset{\curvearrowright}{10}1 \xrightarrow{r(1-\epsilon)} 0011, \quad  
%	1\accentset{\curvearrowright}{10}1 \xrightarrow{r(1-\kappa)} 1011\\[4mm]
%	&0\accentset{\curvearrowleft}{01}0 \xrightarrow{\ell(1+\lambda)} 0100,\quad  
%	1\accentset{\curvearrowleft}{01}0 \xrightarrow{\ell(1-\epsilon)} 1100,\quad  
%	0\accentset{\curvearrowleft}{01}1 \xrightarrow{\ell(1+\epsilon)} 0101, \quad  
%	1\accentset{\curvearrowleft}{01}1 \xrightarrow{\ell(1-\lambda)} 1101
%\end{align*}
%where the curved arrows indicate the jump direction. 
These asymmetric rates furnish us with 
a more generic model that can be studied rigorously and which by virtue of having 
more free parameters has a potentially wider range of applications.

The conventional KLS model is reversible w.r.t. the Ising measure \eqref{Isingperdef}
if and only if there is no bias in the jump rates, i.e., for $r=\ell$. In the generalized
model the situation is slightly more complicated.

\begin{theorem}
\label{Theo:gKLSrev}
The generalized periodic KLS model is reversible w.r.t. the Ising measure \eqref{Isingperdef} if and only if $r=\ell$ and $\kappa=\lambda$.
\end{theorem}

In other words, the generalized periodic KLS model is not reversible, unless it specializes to the conventional symmetric
KLS model where by definition $\kappa=\lambda$.
From a physics perspective the kinetic interaction parameters $\epsilon,\kappa,\lambda$ of the nonreversible generalized KLS model
represent nonreciprocal interactions that are not covered by Newton's third law but are, as effective interactions, ubiquitous in complex
systems, from molecular level in biological cells 
up to collective swarm behavior of living 
beings \cite{Bowi22}. It is generally an important but difficult problem to analyse
the stationary properties of such systems out of thermal equilibrium. Usually such models
are studied using Monte-Carlo simulations which brings about mathematically uncontrolled numerical errors that may
be hard to estimate. Hence it is
useful to have at hand non-reciprocal toy models for which the {\it exact} invariant
measure is known. 
Below we focus on the parameter range $\kappa, \lambda, \epsilon \in (-1,1)$ where the
process is ergodic for any fixed total particle number.
If $|\kappa|=1$, $|\lambda|=1$, or $|\epsilon|=1$, the process may exhibit transient or frozen configurations which we do not consider.

\subsection{Open boundaries}

Periodic boundary conditions are usually considered to study bulk
properties of macroscopically large systems. However, one-dimensional
real systems do not usually consist of an extremely large number of particles.
Moreover, they are frequently attached at their boundaries to external reservoirs
with which particles can be exchanged, examples being
molecular motors moving along a macromolecule such as
DNA or RNA \cite{Chow24}. In this section we establish the existence
of boundary processes that describe such reservoirs and that keep the model 
tractable rigorously
at least along some nontrivial parameter manifold. The strategy is to postulate
a certain form of the invariant measure and of the boundary rates  for particle exchange
processes with free parameters and to determine conditions on these parameters
for the postulated measure to be invariant w.r.t. the postulated boundary processes.
To avoid overly heavy notation
we use the same symbols for the invariant measure and related
quantities as in the periodic case.

A common approach to deal with reservoirs is to allow particles to enter and leave 
the lattice at the
boundary sites $1$ and $L$ respectively as in the open ASEP discussed
in the previous section. In such a setting the translation-invariant periodic 
Ising measure \eqref{Isingperdef} cannot be expected to be invariant. Instead, 
our starting point is the non-translation invariant
Ising measure
\begin{equation}
\pi^\ast_{L}(\bfeta) = 
\frac{1}{Z_{L}}\exp{\left(
\displaystyle - J \sum_{k=1}^{L-1} \eta_k \eta_{k+1} +
\varphi \sum_{k=2}^{L-1} \eta_k + \varphi_{-} \eta_1
 + \varphi_{+}\eta_L
\right)}
\label{Isingfreedef}
\end{equation}
which is obtained from \eqref{Isingperdef} by omitting in the Boltzmann weight \eqref{IsingBW}
 the static interaction term $J\eta_L \eta_{1}$
associated with the lattice bond $(L,1)$ and replacing the chemical potential $\varphi$ appearing in
the boundary term
$\varphi (\eta_L + \eta_{1})$  by boundary chemical potentials $\varphi_{\pm}$, see Appendix \ref{App:Ising} for the transfer matrix treatment
of this measure. Below we express the invariant
measure in terms of the parameters
\begin{equation}
y := \frac{1-\epsilon}{1+\epsilon},
%= \mathrm{e}^{- J} 
\quad
z := \mathrm{e}^{\varphi}, \quad
z_{\pm} := \mathrm{e}^{ \varphi_{\pm}}, \quad q_{\pm} := \frac{z_{\pm}}{z}  \label{xydef} 
\end{equation}
as this simplifies various formulas. In equilibrium physics the quantity
$z$ is called fugacity. Below, the parameters $z_{\pm}$ are referred to as boundary fugacities. We denote the set of strictly positive real numbers
by $R^+$ and assume $y,z,z_{-},z_{+}\in\R^+$ to avoid degenerate
limiting cases.

To obtain an open version of the KLS model one expects that also the boundary jump rates and the rates
of insertion and removal have to account for the two-sided short-range 
kinetic bulk interactions.
In analogy to the interaction range of the bulk hopping rates we postulate
boundary rates that only depend on the boundary sites 1 and $L$ and the sites next to them.

\begin{definition}
\label{Def:gKLSopen}
The generalized open KLS model on the one-dimensional finite integer lattice $\Lambda_L := \{1,\dots,L\}$ is defined by the generator \eqref{genASEP}
with bulk jump rates \eqref{bjr} given by \eqref{rGKLS},  \eqref{lGKLS}
for $2\leq k \leq L-2$, boundary jump rates $w_1(\cdot) = r^{-}(\cdot) + \ell^{-}(\cdot)$,  $w_L(\cdot) = r^{+}(\cdot) + \ell^{+}(\cdot)$ given by
\begin{eqnarray}
r^{-}(\bfeta) & = & \eta_{1} \bar{\eta}_{2} \left[ r^{-}_{1} \bar{\eta}_{3}
+ r^{-}_{2} (1-\epsilon) \eta_{3} \right] 
\label{rm} \\
\ell^{-}(\bfeta) & = & z_{-} z^{-1} \bar{\eta}_{1} \eta_{2} 
 \left[ \ell^{-}_{1} \bar{\eta}_{3} + \ell^{-}_{2} (1+\epsilon) \eta_{3} \right]
\label{lm} \\
r^{+}(\bfeta) & = & z_{+} z^{-1} \eta_{L-1} \bar{\eta}_{L}\left[ r^{+}_{1} \bar{\eta}_{L-2} + r^{+}_{2} (1+\epsilon) \eta_{L-2} \right] 
\label{rp}\\
\ell^{+}(\bfeta) & = & \bar{\eta}_{L-1} \eta_{L} \left[ (\ell^{+}_{1}
\bar{\eta}_{L-2} + \ell^{+}_{2} (1-\epsilon) \eta_{L-2} \right] 
\label{lp}
\end{eqnarray}
for the boundary bonds $(1,2)$ and $(L-1,L)$ respectively,
and insertion and removal of particles with rates
\eqref{boundaryrates} given by
\begin{eqnarray}
\alpha(\bfeta) & = & z_{-} \bar{\eta}_{1} (\alpha_1 \bar{\eta}_2 + \alpha_2 \eta_{2}) 
\label{a12} \\
\gamma(\bfeta) & = & \eta_{1} (\gamma_1 \bar{\eta}_2 + \gamma_2 \eta_{2}) 
\label{g12} \\
\beta(\bfeta) & = &   \eta_{L} (\beta_1 \bar{\eta}_{L-1} + \beta_2 \eta_{L-1})
\label{b12} \\
\delta(\bfeta) & = &  z_{+} \bar{\eta}_{L} (\delta_1 \bar{\eta}_{L-1} + \delta_2 \eta_{L-1}) 
\label{d12}
\end{eqnarray}
for the boundary sites $1$ and $L$ respectively.
\end{definition}

The main result for this open generalized KLS model is the invariance
of the Ising measure \eqref{Isingfreedef} with these boundary processes
for judiciously chosen boundary rates.
For a precise statement and proof of the result it is convenient to introduce the parameters
\begin{eqnarray}
f & := & r  - \ell  \label{fdef} \\
c & := & r \kappa - \ell \lambda \label{cdef}
\end{eqnarray}
for the kinetic bulk interactions,
and the parameters
\begin{eqnarray}
d^{\pm}_{i} & := & r^{\pm}_{i} - \ell^{\pm}_{i}
\label{dpmidef} \\
a_1^{-} & := & \alpha_1   - \gamma_1 \label{c1mdef} \\
a_2^{-} & := & \alpha_2   - \gamma_2 y \label{c2mdef} \\
a_1^{+} & := & \beta_1 - \delta_1 \label{c1pdef} \\
a_2^{+} & := & \beta_2 y - \delta_2  \label{c2pdef} 
\end{eqnarray}
for the kinetic boundary interactions.
For $\epsilon=0$ we also use $c:= r \kappa - \ell \lambda$ and the bias parameter $f:= r-\ell$.
\begin{theorem}
\label{Theo:gKLSopen} (1) Let the interaction parameter $\epsilon$ be in the range $0<|\epsilon|<1$.
For the open generalized KLS model 
the Ising measure \eqref{Isingfreedef} is invariant for any $z\in\R^+$ 
if and only if all of the following conditions (A) - (C) hold:\\
\emph{(A)} The boundary rates
\eqref{rm} - \eqref{d12} are on the parameter manifold defined by
\begin{eqnarray}
d^{\pm}_{1} & = & \frac{(f \epsilon + c) (1+\epsilon) + q_{\pm}^{-1} (f \epsilon - c) (1-\epsilon)}{2 \epsilon} 
\label{dpm1} \\
d^{\pm}_{2} & = & \frac{(f \epsilon + c)  + q_{\pm}^{-1} (f \epsilon - c) }{2 \epsilon}.
\label{dpm2} \\
a_1^{\pm} & = &  \frac{d^{\pm}_{1}}{1+z_{\pm}} 
\label{apm1}\\
a_2^{\pm}  & = & a_1^{\pm}
+ (q_{\pm} - 1) \frac{(f \epsilon + c) (1+\epsilon)}{2 \epsilon z_{\pm}} 
\label{apm2} 
\end{eqnarray}
where $q_{\pm} := z_{\pm} z^{-1}$.\\
\emph{(B)} The boundary fugacities $z_{\pm}$
satisfy the relations
\begin{eqnarray}
0 & = & 2 \epsilon a_2^{-} (z_{-} + y^{-1}) 
- (q_{-}  (f \epsilon + c)  + (f \epsilon - c)) (1+\epsilon)  \label{cH4} \\[2mm]
% % % % % % % % % % % % % % %
0 & = & 2 \epsilon a_2^{+} (z_{+} + y^{-1}) 
- (q_{+}  (f \epsilon + c)  + (f \epsilon - c)) (1+\epsilon).
\label{cH7}
\end{eqnarray}  
\emph{(C)} The compatibility condition
\begin{equation}
a_1^{+} z_{+}  = a_1^{-} z_{-}
\label{cH1}
\end{equation}
between the boundaries holds.\\
(2) For the open generalized KLS model with
$\epsilon=0$ and $c\neq 0$
the Bernoulli product measure 
\begin{equation}
\pi^\ast_{L}(\bfeta) = 
\frac{z_{-}^{\eta_1}}{1+z_{-}} \left(\prod_{k=2}^{L-1} \frac{z^{ \eta_k}}{1+z}\right) \frac{z_{+}^{\eta_L}}{1+z_{+}}
\label{Isingfreedef2}
\end{equation}
 is invariant for any $z\in\R^+$ 
if and only if all of the following conditions hold:\\
\emph{(A0)} The boundary rates
\eqref{rm} - \eqref{d12} of the open generalized KLS model 
with $\epsilon=0$ are on the parameter manifold
defined by
\begin{eqnarray}
d^{\pm}_{1} & = & f + \frac{c}{1+z} 
\label{dpm10} \\
d^{\pm}_{2} & = & f - \frac{c}{1+z^{-1}} 
\label{dpm20} \\
a^{\pm}_1 & = & \frac{d^{\pm}_{1}}{1+z}
\label{apm10} \\
a^{\pm}_2 & = & a^{\pm}_1 - \frac{c}{1+z} 
\label{apm20}
\end{eqnarray}
where $f=r-\ell$, $c=r\kappa - \ell \lambda$.\\
\emph{(B0)}
The boundary fugacities
\begin{eqnarray}
z_{-}=z_{+}=z
\label{cH470}
\end{eqnarray}  
are equal to the bulk fugacity.\\
\emph{(C0)} The compatibility condition
\begin{equation}
a_1^{+}  = a_1^{-}
\label{cH10}
\end{equation}
between the boundaries holds.\\
(3) For the open generalized KLS model with
$\epsilon=0$ and $c= 0$
the Bernoulli product measure \eqref{Isingfreedef2}
 is invariant for any $z,z_{-},z_{+}\in\R^+$ 
if and only if the boundary rates are on the manifold defined by
\begin{eqnarray}
 d^{\pm}_{1} & = &  f 
\frac{1+ z_{\pm}^{-1}}{1 + z^{-1}} 
\label{d1pm00}  \\
 d^{\pm}_{2} & = &  f 
\frac{1+ z_{\pm}^{-1}}{1 + z^{-1}}  \label{d2pm00} \\
a_1^{\pm}  & = &  f \frac{z_{\pm}^{-1}}{1 + z^{-1}}  \label{a1pm00} \\
a_2^{\pm} & = &  f \frac{z^{-1}}{1 + z^{-1}} .
\label{a2pm00}
\end{eqnarray}
and the compatibity condition \eqref{cH1} holds.
\end{theorem}

\begin{remark}
When $\epsilon$ and $c$ are sent to 0 simultaneously, the selected limiting boundary fugacities $z_{\pm}$ (which are free parameters for $\epsilon=c=0$) depend on the way this limit is taken.
\end{remark}
\begin{remark}
For $\epsilon=\kappa=\lambda=0$ the bulk dynamics specializes to the ASEP and the invariant measure is the Bernoulli product measure \eqref{Isingfreedef2} with bulk density $\rho=z/(1+z)$ and boundary densities $\rho_{\pm}
= z_{\pm}/(1+z_{\pm})$. The open ASEP with the \emph{homogeneous} Bernoulli product measure with $z_{-}=z_{+}=z$ as unique invariant measure
is given in terms of the notation of \cite{Derr93a} by $r=p$, $\ell=q$, $r^{\pm}_1 = 
r^{\pm}_2 = p$, $\ell^{\pm}_1 = \ell^{\pm}_2 = q$, $\alpha_1 = 
\alpha_2 = \alpha (1-\rho)/\rho$ 
$\gamma_1 = 
\gamma_2 = \gamma$
$\beta_1 = \beta_2 = \beta$, $\delta_1 = 
\delta_2 = \delta (1-\rho)/\rho$.
The boundary processes in the present
model 
are more general as only the difference relations
$d_1^\pm = d_2^\pm = r-\ell$ and 
$a_1^\pm = a_2^\pm = (r-\ell)(1-\rho)$
need to be satisfied. Consequently the boundary densities
of the Bernoulli product measure may be different from the bulk density $\rho$.
\end{remark}

\begin{theorem}
\label{Theo:gKLSopenrev}
The open generalized KLS model is reversible w.r.t. the Ising measure \eqref{Isingfreedef} if and only if $r=\ell$ and $\kappa=\lambda$ and $d_1^{-}=d_2^{-}=a_1^{-}=a_2^{-}=a_1^{+}=a_2^{+}=0$.
\end{theorem}

Theorem \ref{Theo:gKLSopenrev} implies that the generalized KLS
model is irreversible also for open boundaries, unless it reduces to
the non-driven conventional KLS model.

\subsection{Stationary current and boundary-induced phase transitions}

The main quantity that expresses the nonequilibrium nature of the model is
the stationary current \eqref{jLdef}.
The computation of this expectation is straightforward using the well-known
properties of the Ising measure summarized in Appendix \ref{App:Ising}. Of particular importance are the eigenvalues
$\lambda_{0}$ and $\lambda_{1}$ defined \eqref{eigenvalues}
of the
transfer matrix \eqref{transfermatrix} that encodes the
invariant measure of the KLS model and the auxiliary quantities
$Y_n$ defined in \eqref{Ydef}.
Using translation invariance then yields
in terms of the bulk parameters \eqref{fdef} and \eqref{cdef} the exact expression
\bea 
j(L) 
& = & (f + c) \exval{\eta_k \bar{\eta}_{k+1}}_L +  (f \epsilon - c) \exval{\eta_{k-1}\eta_{k} \bar{\eta}_{k+1}}_L 
- (f \epsilon + c) \exval{\eta_{k-1} \bar{\eta}_{k} \eta_{k+1}}_L   \\
& = & \left[ f + c +  (f (1-\epsilon) - 2 c)z^2 y \frac{Y_{L-1}}{Y_L} 
 - (f \epsilon + c + (f \epsilon - c) z^2 y^2) x^2 
\frac{ Y_{L-2}}{Y_L} \right] \rho(L)  
\label{genKLScurrent}
\eea
for a finite periodic system where in the second equality the expressions
\eqref{GIsing3} - \eqref{HIsing32} are inserted. In the thermodynamic limit $L\to\infty$ this expression
reduces to
\begin{eqnarray}
j & := & \lim_{L\to\infty} j(L) \\
& = & \left[  (f + c) +  ((f \epsilon - c) - (f + c))z^2 y \lambda_{0}^{-1} 
 - ((f \epsilon + c) + (f \epsilon - c) z^2 y^2) z^2 
\lambda_{0}^{-2} \right] \rho 
\end{eqnarray}

\begin{figure}[H]
\centering
\includegraphics[width=8cm]{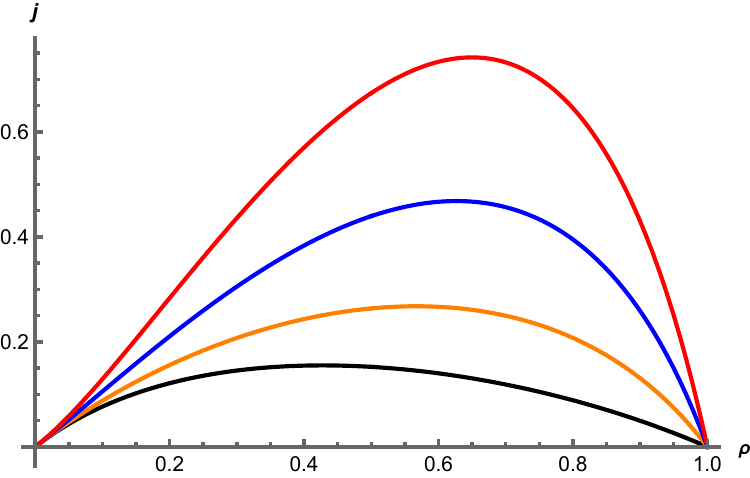}
\caption{Stationary current $j $  as a function of the particle density 
$\rho $ for repulsive static interaction with $\epsilon = 2/3$ and with kinetic interaction parameters $\ell =0.9,  \lambda = 1/9$ 
and different values of $r$ and $\kappa$. Curves from the bottom to top: 
$(r,\kappa) = (1.606, 0.650), (1.720, 0.389), (1.900, 1/9), (2.140, -0.123)$.}
	\label{current_y_05}
\end{figure}

\begin{figure}[H]
	\centering
	\includegraphics[width=8cm]{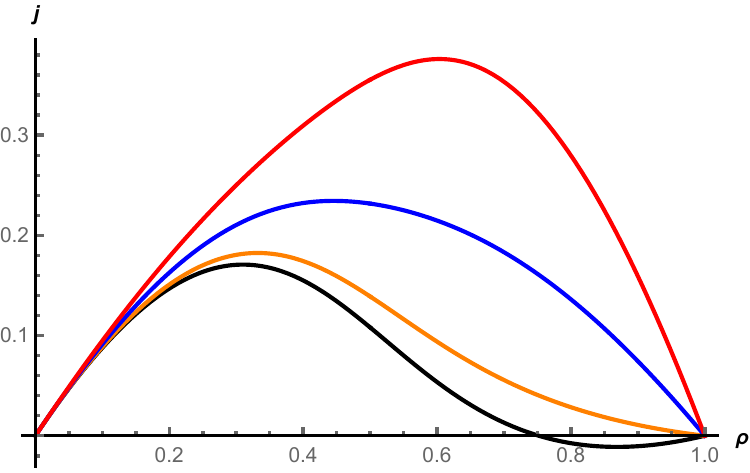}
\caption{Stationary current $j $  as a function of the particle density 
$\rho $ for attractive static interaction with $\epsilon = - 2/3$ and with kinetic interaction parameters $\ell =4.5,  \lambda = -7/9$ 
and different values of $r$ and $\kappa$. Curves from the bottom to top: 
$(r,\kappa) = (4.3, - 0.697), (4.6, -0.722), (5.5, -0.778), (6.7, -0.825)$.}
	\label{current_y_5}
\end{figure}

\begin{figure}[H]
	\centering
	\includegraphics[width=8cm]{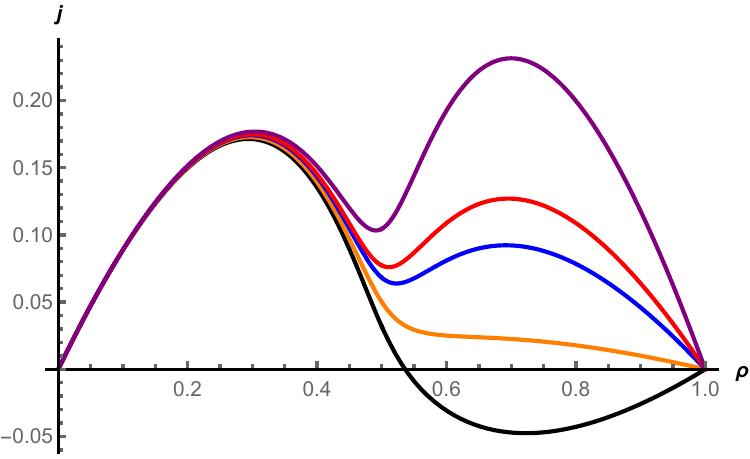}
\caption{Stationary current $j$  as a function of the particle density $\rho$ for 
strong attractive static interaction with $\epsilon = - 0.980$ with kinetic interaction parameters $\ell =150.75$, 
$\lambda = -0.987$ 
and different values of $r$ and $\kappa$. Curves from the bottom to top: 
$(r,\kappa) = (56.55, - 0.982 ), (66.65, -0.985), (76.75, -0.987), (81.8, -0.988), (96.95, -0.989)$.}
	\label{current_y_100}
\end{figure}

By the general properties of the Ising measure
the stationary current of the open system deviates from the periodic case only
by terms that are exponentially small in system size with the same parameter
$\ln{(\lambda_{1}/\lambda_{0})}$ that determines the magnitude of the finite-size
corrections of the finite periodic system relative to the infinite system.
To illustrate the result the stationary current in the thermodynamic limit
is plotted in Figs. \ref{current_y_05} - 
\ref{current_y_100} as a function of the particle density for some choices of the rates.

%\begin{figure}[H]
%	\centering
%	\includegraphics[width=8cm]{Fig_Appen/y=0.5_d_001_02_05_09.pdf}
%	\caption{Stationary current $j$ with $\dfrac{1+r^{\yesp 
%\star}}{1+r^{\star\nop\yesp}} = 5 \text{ (repulsive
%			interaction)}, r^\star = 2, \ell^\star =1,  \ell^{\star\yesp} = 0.5$ and 
%different values of $r^{\yesp\star}$ as a function of the particle density 
%$\rho$. 
%Curves from the bottom to top: $r^{\yesp\star} = 0.01, 0.2, 0.5, 0.9$.}
%	\label{current_y_05}
%\end{figure}
%

%\begin{figure}[H]
%	\centering
%	\includegraphics[width=8cm]{Fig_Appen/y=5}
%	\caption{Stationary current $j$ with $\dfrac{1+r^{\yesp \star}}{1+r^{\star\nop\yesp}} = 0.2$ 
%(attractive interaction), $r^\star = 2, \ell^\star =1,  \ell^{\star\yesp} = 0.5$ and different values of $r^{\yesp\star}$ as a function of the particle density $\rho$. Curves from the bottom to top: $r^{\yesp\star} = 0.1, 0.2, 0.5, 0.9$.}
%	\label{current_y_5}
%\end{figure}

%\begin{figure}[H]
%	\centering
%	\includegraphics[width=8cm]{Fig_Appen/y=100}
%	\caption{Stationary current $j$ with $\epsilon = - 0.98$
%%$\dfrac{1+r^{\yesp \star}}{1+r^{\star\nop\yesp}} = 0.01$ 
%(strong attractive interaction), $r^\star = 2, \ell^\star =1,  \ell^{\star\yesp} = 0.5$ and different values of $r^{\yesp\star}$ as a function of the particle density $\rho$. Curves from the bottom to top: $r^{\yesp\star} = 0.1, 0.3, 0.5, 0.6, 0.9$. }
%	\label{current_y_100}
%\end{figure}
%	\begin{equation}
%		\begin{cases}
%			r^\star &=\ r(1+\kappa) \\
%			r^{\yesp \star} &=\ \dfrac{1+\epsilon}{1+\kappa}-1 \\
%			r^{\star\nop\yesp} &=\ \dfrac{1-\epsilon}{1+\kappa}-1  \\
%			\ell^\star &=\ \ell(1+\lambda) \\
%			\ell^{\star\yesp} &=\ \dfrac{1+\epsilon}{1+\lambda}-1 \\
%			\ell^{\yesp\nop\star} &=\ \dfrac{1-\epsilon}{1+\lambda}-1.  
%		\end{cases}
%	\end{equation}
%
%\be 
% \frac{1+ r^{\yesp \star}}{1+ r^{\star\nop\yesp}} = \dfrac{1+\epsilon}{1-\epsilon}
%\ee
%
%\bea
%r & = & r^\star + \frac{r^{\yesp \star}+r^{\star\nop\yesp}}{2} \\
%\ell & = & \ell^\star + \frac{\ell^{\star\yesp }+\ell^{\yesp\nop\star}}{2} \\
% \kappa & = & - \frac{r^{\yesp \star}+r^{\star\nop\yesp}}{2+ r^{\yesp \star}+r^{\star\nop\yesp}} \\
%\lambda & = & - \frac{\ell^{\star\yesp }+\ell^{\yesp\nop\star}}{2+ \ell^{\star\yesp }+\ell^{\yesp\nop\star}} \\
%%1+ \kappa & = & \frac{2}{2+ r^{\yesp \star}+r^{\star\nop\yesp}} \\
%\epsilon & = &  
%\frac{r^{\yesp \star} - r^{\star\nop\yesp}}{2+ r^{\yesp \star}+r^{\star\nop\yesp}} 
%\, = \, \frac{\ell^{\star\yesp }-\ell^{\yesp\nop\star}}{2+ \ell^{\star\yesp }+\ell^{\yesp\nop\star}} 
%\eea
%

\begin{remark}
The stationary current may exhibit current reversal, i.e., a change sign as a function of density, as shown in Figs. \ref{current_y_5} and \ref{current_y_100}. Analytically, 
this is most easily seen for $\epsilon=0$ (no static interaction) 
where %\bea
%\langle r_k \rangle_{\rho} 
%& =	& r \rho (1-\rho) [1+\kappa (1 -  2 \rho) ] \\
%\langle \ell_{k+1}  \rangle_{\rho} 
%& =	& \ell \rho (1-\rho) [1+\lambda (1 -  2 \rho) ] \\
%\eea
\be 
j =\rho (1-\rho) \left[f  + c (1-2\rho)\right].
\ee
For fixed transition rates the current changes sign at density
\begin{equation}
\rho_0 = \frac{1}{2} \frac{f+c}{c}
%\in (0,1) 
\end{equation}
so that current reversal may occur for parameters $f$ and $c$ taken
such that $-1 < f/c < 1$. This is in contrast to the conventional KLS model
where the sign of the stationary current is determined for all densities by the sign of the
hopping bias $f$ alone. This can be seen by noting that for $\kappa=\lambda$ the stationary current \eqref{jLdef} is given by $j = f j^{+}$
where $j^{+} = \exval{r_k}$ represents the stationary current from the site 
$k$ to site $k+1$ which is positive
by definition. Current reversal has been observed previously
in a periodic zero-range type processes without exclusion \cite{Chat17} and
in a family of multi-species exclusion processes \cite{Chat18}.
\end{remark}

\begin{figure}[H]
	\centering
	\includegraphics[width=8cm]{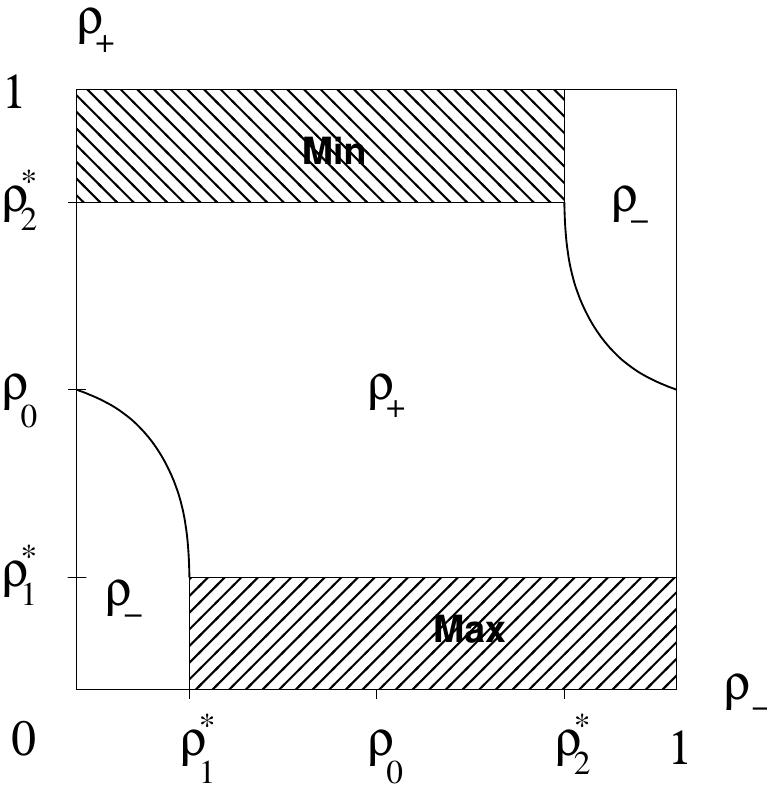}
\caption{Schematic phase diagram of the generalized KLS model for current reversal, with left boundary density $\rho_{-}$ on the $x$-axis and right boundary density $\rho_{+}$ on the $y$-axis. The minimal (maximal) current phase is marked by Min (Max). Both are bounded by straight second-order phase transition lines.
In the central region marked by $\rho_{+}$ the bulk density is equal to the right boundary
density. Curved first-order phase transition lines separate it from two regions with bulk density equal to the left boundary density
(lower left and upper right regions marked by $\rho_{-}$).}
	\label{phasediagram}
\end{figure}

\begin{remark}
In the presence of current reversal with two densities $\rho_1^\ast$ where $j'(\rho)=0$ the extremal current principle \eqref{ex} predicts the universal phase diagram
shown in Fig. \ref{phasediagram}. There are two first-order phase transition lines 
determined by the condition $j(\rho_{-})= j(\rho_{+})$ and two second-order 
phase transition lines, one bounding the maximal current phase
with $\rho = \rho_1^\ast$ in the range $\rho_1^\ast \leq \rho_{-} \leq 1$,
$0 \leq \rho_{+} \leq \rho_1^\ast$ and the other one 
bounding the minimal current phase
with $\rho = \rho_2^\ast$ in the range $0 \leq \rho_{-} \leq \rho_2^\ast$,
$\rho_2^\ast \leq \rho_{+} \leq 1$.  The open ASEP only exhibits the phases that appear in the lower left square $[0,\rho_0]^2$.
\end{remark}

\begin{remark}
In the open system the stationary current is negative for $\rho_{+}>\rho_{0}$
independently of $\rho_{-}$. Thus the current flows uphill against the external density gradient in the whole region defined by $\rho_{-}>\rho_{+}>\rho_{0}$.
Likewise, the current flows uphill against the external density gradient in the whole region defined by $\rho_{-}<\rho_{+}<\rho_{0}$.
\end{remark}

\subsection{Proofs}

For the proofs it is convenient to express the generator in terms of the intensity matrix \cite{Ligg10}.
To this end, we introduce the two-dimensional matrices
\begin{equation}
\sigma^{+} := \left(\ba{cc} 0 & 1 \\ 0 & 0 \ea\right), \quad
\sigma^{-} := \left(\ba{cc} 0 & 1 \\ 0 & 0 \ea\right), \quad
\hat{n} := \left(\ba{cc} 1 & 0 \\ 0 & 0 \ea\right), \quad
\hat{v} := \left(\ba{cc} 0 & 1 \\ 0 & 0 \ea\right) 
\end{equation}
and the two-dimensional unit matrix $\mathds{1}$. From these we construct the
local operators $a_k := \mathds{1}^{(k-1)} \otimes a \otimes \mathds{1}^{(L-k)}$
via the Kronecker product $\otimes$.

We also define , the
$L$-fold Kronecker product $\langle s | := (1,1)^{\otimes L}$, the $2$-dimensional row vector $(1,1)$ and the
column vector $| s \rangle = \langle s |^T$ obtained from
$\langle s |$ by transposition. This is the column
vector of dimension $2^L$ with all it s components equal to 1
which satisfies the relations
\begin{equation}
\sigma_{k}^+ | s \rangle = \hat{v}_{k} | s \rangle , \quad
\sigma_{k}^- | s \rangle = \hat{n}_{k} | s \rangle, \quad | \mu \rangle = \hat{\mu} | s \rangle ,
\label{skmmu}
\end{equation}
see \cite{Schu01,Schu19} for details.

Next, with the local operators we define
the diagonal bulk and boundary jump rate matrices
\bea
\hat{r}_{k} & = & r \hat{n}_{k} \hat{v}_{k+1} \left[ 
(1+\kappa) (\hat{v}_{k-1} - \hat{n}_{k+2}) + (1+\epsilon) \hat{n}_{k-1}
+ (1-\epsilon) \hat{n}_{k+2}\right] , \quad
2 \leq k \leq L-2 
\label{hrkdef} \\
\hat{\ell}_{k+1} & = & \ell \hat{v}_{k} \hat{n}_{k+1} \left[
(1+\lambda)(\hat{v}_{k-1} - \hat{n}_{k+2}) +
(1 - \epsilon) \hat{n}_{k-1} + (1+\epsilon) \hat{n}_{k+2} 
\right] , \quad
2 \leq k \leq L-2 
\label{hlkdef} \\
\hat{r}^{-} & = & \hat{n}_{1} \hat{v}_{2}
(r^{-}_{1} \hat{v}_3 + r^{-}_{2} (1-\epsilon) \hat{n}_3)
\label{hrmdef} \\
\hat{\ell}^{-} & = & z_{-} z^{-1} \hat{v}_{1} \hat{n}_{2} 
(\ell^{-}_{1} \hat{n}_3 + \ell^{-}_{2} (1+\epsilon) \hat{n}_3)
\label{hlmdef}  \\
\hat{r}^{+} & = & z_{+} z^{-1} \hat{v}_{L-1} \hat{n}_{L} (r^{+}_{1} \hat{v}_3 + r^{+}_{2} (1+\epsilon) \hat{n}_3)
\label{hrpdef}  \\
\hat{\ell}^{+} & = & \hat{v}_{L-1} \hat{n}_{L} (\ell^{+}_{1} \hat{n}_3 + \ell^{+}_{2} (1-\epsilon) \hat{n}_3) 
\label{hlpdef} 
\eea
and the diagonal boundary rate matrices
\bea
\hat{\alpha} & = & z_{-} \hat{v}_{1} (\alpha_1 \hat{v}_{2} + \alpha_2 \hat{n}_{2}) \\
\hat{\gamma} & = & \hat{n}_{1} (\gamma_1 \hat{v}_{2} + \gamma_2 \hat{n}_{2})  \\
\hat{\beta} & = &   \hat{n}_{L} (\beta_1 \hat{v}_{L-1} + \beta_2 \hat{n}_{L-1}) \\
\hat{\delta} & = & z_{+} \hat{v}_{L} (\delta_1 \hat{v}_{L-1} + \delta_2 \hat{n}_{L-1}) .
\label{Hdef}
\eea
The intensity matrix representing the generator \eqref{Def:gKLSopen} is then 
in the quantum Hamiltonian formulation \cite{Alex78,Sudb95,Schu01} given by
\be 
H = - \left(b^{-} + h^{-} + \sum_{k=2}^{L-2} h_{k} + h^{+} + b^{+}\right)
\label{HKLSopen}
\ee
where
\begin{eqnarray}
h_{k} 
& = & \left( \sigma_{k}^{+} \sigma_{k+1}^{-} 
- \hat{n}_{k} \hat{v}_{k+1} \right) \hat{r}_{k} 
+ \left(\sigma_{k}^{-} \sigma_{k+1}^{+} 
- \hat{v}_{k} \hat{n}_{k+1} \right) \hat{\ell}_{k+1} , \quad
2 \leq k \leq L-2 
\label{hkdef} \\
h^{-} 
& = & \left( \sigma_{1}^{+} \sigma_{2}^{-} 
- \hat{n}_{1} \hat{v}_{2} \right) \hat{r}^{-}
+ \left( \sigma_{1}^{-} \sigma_{2}^{+} 
- \hat{v}_{1} \hat{n}_{2} \right) \hat{\ell}^{-}
\label{hmdef} \\
h^{+} 
& = & \left( \sigma_{L-1}^{+} 
\sigma_{L}^{-} 
- \hat{n}_{L-1} \hat{v}_{L} \right) \hat{r}^{+}
+ \left(\sigma_{L-1}^{-} 
\sigma_{L}^{+} 
- \hat{v}_{L-1} \hat{n}_{L} \right) \hat{\ell}^{+}
\label{hpdef} \\
b^{-} & = & (\sigma^{-}_1 - \hat{v}_{1}) \hat{\alpha} + 
(\sigma^{+}_1 - \hat{n}_{1}) \hat{\gamma} 
\label{bmdef} \\
b^{+} & = & (\sigma^{-}_{L} - \hat{v}_{L}) \hat{\delta} + 
(\sigma^{+}_{L} - \hat{n}_{L}) \hat{\beta}  .
\label{bpdef} 
\end{eqnarray}
For periodic boundary conditions the quantum
Hamiltonian form of the generator is simply given by
\be 
H = -  \sum_{k=1}^{L} h_{k} 
\label{HKLSper}
\ee
with lattice indices understood to be modulo $L$.

\subsubsection{Proof of Theorem \ref{Theo:gKLSrev}}

The central idea is to use the ground state
transformation 
\begin{equation}
\tilde{a}_k := \hat{\mu}^{-1} a_k \hat{\mu}  
\label{gstdef}
\end{equation}
of local operators where
\begin{equation}
\hat{\mu} = \prod_{k=1}^{L}
\left[1 + \left(y-1\right) \hat{n}_k 
\hat{n}_{k+1}\right] z^{\hat{n}_k} 
\label{hmudef}
\end{equation}
is the diagonal invariant measure matrix 
with inverse 
\begin{equation}
\hat{\mu}^{-1} = \prod_{k=1}^{L}
\left[1 + \left(y^{-1}-1\right) \hat{n}_k 
\hat{n}_{k+1}\right] z^{-\hat{n}_k} .
\label{hmuinv}
\end{equation}
Reversibility -- which is the condition of detailed balance --
then reads \cite{Schu19}
\begin{equation}
H = \tilde{H}^T
\label{detbalH}
\end{equation}
where $\tilde{H} := \hat{\mu}^{-1} H \hat{\mu}$ is the ground state
transformation of $H$ given by \eqref{HKLSper}.
Since detailed balance must hold for all pairs $\bfeta,\bfeta'$
of configurations and all transitions between those states it is necessary and sufficient
that \eqref{detbalH} holds for every local term appearing in \eqref{Hdef}. 

\paragraph{(i) Ground state transformation:}

\begin{proposition}
\label{Prop:gst1}
For the local operators $\sigma^{\pm}_{k}$ the ground state transformation with the diagonal measure matrix 
\eqref{hmudef} yields
\bea 
\hat{\mu}^{-1} \sigma_k^{\pm} \hat{\mu}  
& = & z^{\pm 1} \sigma_k^{\pm} y^{\pm(\hat{n}_{k-1}+\hat{n}_{k+1})}  , \quad 1 \leq k \leq L
\eea
with $k\in\{1,\dots,L\}$ taken modulo $L$.
\end{proposition}

The proof is elementary and follows from the commutation properties
$a_k b_l - b_l a_k = 0$ for $k\neq l$ , $a_k b_k - b_k a_k = (a b - b a)_k$
of the Kronecker product 
for any pair of two-dimensional matrices $a,b$
and the matrix multiplication relations 
\be 
%y^{\hat{n}} = \frac{1+\epsilon - 2\epsilon \hat{n}}{1+\epsilon}, \quad
\hat{n} \sigma^{-} = \sigma^{-} \hat{v} = \sigma^{-}, \quad
\hat{v} \sigma^{+} = \sigma^{+} \hat{n} = \sigma^{+}, \quad
\hat{n} y^{\hat{n}} = y \hat{n} .
%, \quad \epsilon = \frac{1-y}{1+y} 
\ee
With these preparations Theorem \ref{Theo:gKLSrev} can be proved as outlined above.

\paragraph{(ii) Proof of reversibility:} From Proposition \ref{Prop:gst1} the ground state transformation \eqref{gstdef} yields
\begin{eqnarray}
\hat{\mu}^{-1} \sigma_{k}^{+} \sigma_{k+1}^{-} \hat{\mu}
& = &  \sigma_k^{+} \sigma_{k+1}^{-} y^{ \hat{n}_{k-1} - 
\hat{n}_{k+2}} \\
%%%%%%%%%%%%%%%%%%%%%%%%%%%%%%%%%
\hat{\mu}^{-1} \sigma_{k}^{-} \sigma_{k+1}^{+} \hat{\mu}
& = &  \sigma_k^{-} \sigma_{k+1}^{+} y^{ \hat{n}_{k+2} - 
\hat{n}_{k-1}} 
\end{eqnarray}
and the transposition property $(\sigma^{\pm})^T= \sigma^{\mp}$
then leads to
\bea
\tilde{h}_{k} 
& = & 
r \sigma_{k}^{-} \sigma_{k+1}^{+} y^{ \hat{n}_{k-1} - \hat{n}_{k+2}} \left[ 
(1+\kappa) (\hat{v}_{k-1} - \hat{n}_{k+2}) + (1+\epsilon) \hat{n}_{k-1}
+ (1-\epsilon)  \hat{n}_{k+2}\right] \nonumber \\ 
& & - r \hat{n}_{k} \hat{v}_{k+1} \left[ 
(1+\kappa) (\hat{v}_{k-1} - \hat{n}_{k+2}) + (1+\epsilon) \hat{n}_{k-1}
+ (1-\epsilon)  \hat{n}_{k+2}\right] \nonumber \\ 
& & + \ell  \sigma_{k}^{+} \sigma_{k+1}^{-} y^{\hat{n}_{k+2} - \hat{n}_{k-1}}
 \left[
(1+\lambda)(\hat{v}_{k-1} - \hat{n}_{k+2}) +
(1 + \epsilon) \hat{n}_{k-1} + (1 - \epsilon) \hat{n}_{k+2} 
\right]  \nonumber \\ 
& & - \ell \hat{v}_{k} \hat{n}_{k+1} \left[
(1+\lambda)(\hat{v}_{k-1} - \hat{n}_{k+2}) +
(1 + \epsilon) \hat{n}_{k-1} + (1 - \epsilon) \hat{n}_{k+2} 
\right]  \nonumber \\
& = & 
r \sigma_{k}^{-} \sigma_{k+1}^{+}  \left[ 
(1+\kappa) (\hat{v}_{k-1} - \hat{n}_{k+2}) + (1-\epsilon) \hat{n}_{k-1}
+ (1+\epsilon)  \hat{n}_{k+2}\right] \nonumber \\ 
& & - \ell \hat{v}_{k} \hat{n}_{k+1} \left[
(1+\lambda)(\hat{v}_{k-1} - \hat{n}_{k+2}) +
(1 + \epsilon) \hat{n}_{k-1} + (1 - \epsilon) \hat{n}_{k+2} 
\right] \nonumber \\
& & + \ell  \sigma_{k}^{+} \sigma_{k+1}^{-}
 \left[
(1+\lambda)(\hat{v}_{k-1} - \hat{n}_{k+2}) +
(1 - \epsilon) \hat{n}_{k-1} + (1 + \epsilon) \hat{n}_{k+2} 
\right]  \nonumber \\ 
& & - r \hat{n}_{k} \hat{v}_{k+1} \left[ 
(1+\kappa) (\hat{v}_{k-1} - \hat{n}_{k+2}) + (1+\epsilon) \hat{n}_{k-1}
+ (1-\epsilon)  \hat{n}_{k+2}\right] .
 \label{hkrev}
\eea
Comparing term by term with \eqref{hmdef} one sees that $\tilde{h}_{k} 
= h_k$ is satisfied if and only if $r=\ell$ and $\kappa=\lambda$. \qed

\subsubsection{Proof of Theorem \ref{Theo:gKLSopen}}

For finite state space invariance of a measure is equivalent to the eigenvalue
condition 
\begin{equation}
H | \mu \rangle =0
\label{Hmuinv}
\end{equation}  
for the quantum Hamiltonian representation of the intensity matrix. 
The proof of this property involves two steps: (i) Compute the action of the local intensity matrices on the Ising measure \eqref{Isingfreedef} by using the ground state transformation and observing the relations
\eqref{skmmu} (see Lemma \ref{Lmm:Hmu} below).
(ii) Fix parameters such that cancellation of ``unwanted'' terms (see below) proves invariance by establishing the
conditions for invariance
of the modified Ising measure \eqref{Isingfreedef}
which yields the diagonal invariant measure matrix
\begin{equation}
\hat{\mu} = z_{-}^{\hat{n}_{1}}
 \left(\prod_{k=2}^{L-1} z^{\hat{n}_{k}} \right) z_{+}^{\hat{n}_{L}}
\prod_{k=1}^{L-1}
\left[1 + \left(y-1\right) \hat{n}_k 
\hat{n}_{k+1}\right] .
\label{hmudef2}
\end{equation}

\paragraph{(i) Ground state transformation:}

\begin{proposition}
\label{Prop:gst2}
For the local operators $\sigma^{\pm}_{k}$ the ground state transformation with the diagonal measure matrix 
\eqref{hmudef2} yields
\bea 
\hat{\mu}^{-1} \sigma_1^{\pm} \hat{\mu}  
& = & z_{-}^{\pm 1} \sigma_1^{\pm} y^{\pm \hat{n}_{2}} 
\label{mus1}  \\
%%%%%%%%%%%%%%%%%%%%%%%%%%%%%%%%%
%\hat{\mu}^{-1} \sigma_1^{-} \hat{\mu}  
%& = & z^{-1} \sigma_1^{-} y^{- \hat{n}_{2}} \nonumber \\
%%%%%%%%%%%%%%%%%%%%%%%%%%%%%%%%%%
\hat{\mu}^{-1} \sigma_k^{\pm} \hat{\mu}  
& = & z^{\pm 1} \sigma_k^{\pm} y^{\pm(\hat{n}_{k-1}+\hat{n}_{k+1})}  , \quad 2 \leq k \leq L-1
 \\
%%%%%%%%%%%%%%%%%%%%%%%%%%%%%%%%%
%\hat{\mu}^{-1} \sigma_k^{-} \hat{\mu}  
%& = & z^{-1} \sigma_k^{-} y^{ - \hat{n}_{k-1} - \hat{n}_{k+1}}  , \quad 2 \leq k \leq L-1
%\nonumber \\
\hat{\mu}^{-1} \sigma_L^{\pm} \hat{\mu}  
& = & z_{+}^{\pm 1} \sigma_L^{\pm} y^{\pm \hat{n}_{L-1}} 
\label{musL} .
%\nonumber \\
%%%%%%%%%%%%%%%%%%%%%%%%%%%%%%%%%%
%\hat{\mu}^{-1} \sigma_L^{-} \hat{\mu}  
%& = & x^{-1} \sigma_L^{-} y^{- \hat{n}_{L-1}} \nonumber
\eea
\end{proposition}

The proof requires the same elementary properties of the Kronecker
product as in the proof of Proposition \ref{Prop:gst1}.

\begin{corollary}
The jump terms $\sigma_{k}^{\pm} \sigma_{k+1}^{\mp}$ satisfy the
transformation properties
\begin{eqnarray}
\hat{\mu}^{-1} \sigma_1^{+} \sigma_{2}^{-} \hat{\mu}
& = & \sigma_1^{+} \sigma_2^{-} y^{- \hat{n}_{3}} z_{-} z^{-1}
\label{musma} \\
%%%%%%%%%%%%%%%%%%%%%%%%%%%%%%%%%
\hat{\mu}^{-1} \sigma_1^{-} \sigma_{2}^{+} \hat{\mu}
& = & \sigma_1^{-} \sigma_2^{+} y^{\hat{n}_{3}} z_{-}^{-1} z
\label{musmb} \\
%%%%%%%%%%%%%%%%%%%%%%%%%%%%%%%%%
\hat{\mu}^{-1} \sigma_{k}^{+} \sigma_{k+1}^{-} \hat{\mu}
& = &  \sigma_k^{+} \sigma_{k+1}^{-} y^{ \hat{n}_{k-1} - 
\hat{n}_{k+2}} , \quad 2 \leq k \leq L-2
\label{muska} \\
%%%%%%%%%%%%%%%%%%%%%%%%%%%%%%%%%
\hat{\mu}^{-1} \sigma_{k}^{-} \sigma_{k+1}^{+} \hat{\mu}
& = &  \sigma_k^{-} \sigma_{k+1}^{+} y^{ \hat{n}_{k+2} - 
\hat{n}_{k-1}}  , \quad 2 \leq k \leq L-2
\label{muskb}  \\
\hat{\mu}^{-1} \sigma_{L-1}^{+} \sigma_{L}^{-} \hat{\mu}
& = & \sigma_{L-1}^{+} \sigma_L^{-} y^{\hat{n}_{L-2}} z_{+}^{-1} z
\label{muspa} \\
\hat{\mu}^{-1} \sigma_{L-1}^{-} \sigma_{L}^{+} \hat{\mu}
& = & \sigma_{L-1}^{-} \sigma_L^{+} y^{-\hat{n}_{L-2}} z_{+} z^{-1}.
\label{muspb} 
\end{eqnarray}
\end{corollary}

\paragraph{(ii) A technical lemma:}

\begin{lemma}
\label{Lmm:Hmu}
With $q_{\pm} := z_{\pm} z^{-1}$ and the parameters \eqref{fdef} and \eqref{cdef}
the relations
\begin{eqnarray}
\sum_{k=2}^{L-2} h_{k} |\mu \rangle
& = &  - \left[(f + c) \hat{n}_{2} + (f \epsilon + c) \hat{n}_{1} \hat{v}_{2} \hat{n}_{3} 
+ (f \epsilon - c) \hat{n}_{2} (1 - \hat{v}_{1} \hat{v}_{3}) \right] |\mu \rangle \nonumber \\[-2mm]
& & + \left[ (f + c) \hat{n}_{L-1} 
 + (f \epsilon + c) \hat{n}_{L-2} \hat{v}_{L-1} \hat{n}_{L}
+ (f \epsilon - c) \hat{n}_{L-1} (1
 - \hat{v}_{L-2} \hat{v}_{L}) \right] |\mu \rangle ,
\label{bulksum} 
\end{eqnarray}

\begin{eqnarray}
h^{-} | \mu \rangle
%& = &  r \left[(1 + \kappa) \left( \hat{n}_{2} - \hat{n}_{1} \right) 
%+ (\epsilon + \kappa) \hat{n}_{1} \hat{n}_{3} 
%+ (\epsilon - \kappa) \hat{n}_{2} \hat{n}_{3} 
% - 2 \epsilon \, \hat{n}_{1} \hat{n}_{2} \hat{n}_{3} \right] |\mu \rangle \nonumber \\
%& & 
%- \ell \left[(1 + \lambda) (\hat{n}_{2} - \hat{n}_{1}) 
%+ (\epsilon + \lambda) \hat{n}_{1} \hat{n}_{3} 
%+ (\epsilon - \lambda) \hat{n}_{2} \hat{n}_{3} 
%- 2 \epsilon \, \hat{n}_{1} \hat{n}_{2} \hat{n}_{3} \right] |\mu \rangle 
%\label{lb} \\
& = &  - d^{-}_{1} \hat{n}_{1}  \,  |\mu\rangle
+ q_{-} d^{-}_{1} \hat{n}_{2} \,  |\mu\rangle - (q_{-}-1) d^{-}_{1}\hat{n}_{1}\hat{n}_{2}   \,  |\mu\rangle \nonumber \\
& & + [d^{-}_{1} - d^{-}_{2}(1-\epsilon)] \hat{n}_{1} \hat{n}_{3} \,  |\mu\rangle
 + q_{-} [ - d^{-}_{1} + d^{-}_{2}(1+\epsilon)]  \hat{n}_{2}\hat{n}_{3} \,  |\mu\rangle 
\nonumber \\
& & - [ d^{-}_{1}(1-q_{-})   - d^{-}_{2} (1-\epsilon)
+ d^{-}_{2} q_{-} (1+\epsilon)] \hat{n}_{1}\hat{n}_{2}\hat{n}_{3} \,  |\mu\rangle 
\label{lb} 
\end{eqnarray}

\begin{eqnarray}
h^{+} | \mu \rangle
& = & d^{+}_{1} \hat{n}_{L}
- q_{+} d^{+}_{1} \hat{n}_{L-1}
 + (q_{+} - 1) d^{+}_{1} \hat{n}_{L-1}\hat{n}_{L}
\nonumber \\
& & - [d^{+}_{1}
- d^{+}_{2} (1-\epsilon) ] \hat{n}_{L-2}\hat{n}_{L} - q_{+} [- d^{+}_{1}
+ d^{+}_{2} (1+\epsilon)] \hat{n}_{L-2}\hat{n}_{L-1} 
\nonumber \\
& & - [(q_{+} - 1) d^{+}_{1} + (1-\epsilon) d^{+}_{2} - q_{+} (1+\epsilon) d^{+}_{2}] \hat{n}_{L-2}\hat{n}_{L-1}\hat{n}_{L} 
\label{rb} 
\end{eqnarray}
and
\begin{eqnarray}
b^{-} | \mu \rangle 
%& = &  - (\alpha_1 z_{-} - \gamma_1 z_{-}) \, | s \rangle \\
%& & + (\alpha_1 z_{-} - \gamma_1 z_{-}) (1+z_{-}^{-1}) \hat{n}_{1}  \, | s \rangle \\
%& & +[\alpha_1 z_{-} - \gamma_1 z_{-} - (\alpha_2 z_{-} -  \gamma_2 z_{-} y) ] \hat{n}_{2} \\
%& & - \left[(\alpha_1 z_{-} - \gamma_1 z_{-}) (1+z_{-}^{-1})- (\alpha_2 z_{-} -  \gamma_2 z_{-} y)  (1+y^{-1} z_{-}^{-1})\right] \hat{n}_{1}  
% \hat{n}_{2}\, | s \rangle \\
& = & \left[ - a_1^{-} z_{-}
+ a_1^{-} (1+z_{-}) \hat{n}_1 + (a_1^{-} - a_2^{-}) z_{-} \hat{n}_{2} 
- [a_1^{-}  (1+z_{-}) - a_2^{-} 
(z_{-} + y^{-1})] \hat{n}_1 \hat{n}_{2} \right] | \mu \rangle 
\label{lres} \\
b^{+} | \mu \rangle 
& = & \left[ a_1^{+} z_{+} -  (a_1^{+} - a_2^{+}) z_{+}\hat{n}_{L-1} 
- a_1^{+} (1+z_{+}) \hat{n}_{L} 
+ [a_1^{+} (1+z_{+})  - a_2^{+} (z_{+} + y^{-1})] \hat{n}_{L-1} \hat{n}_{L} \right] | \mu \rangle .
\label{rres}
\end{eqnarray}
hold for all $\epsilon,\kappa,\lambda\in[-1,1]$. 
%For 
%$\kappa,\lambda\in[-1,1]$, $z_{-}=z_{+}=z$ and
%$\epsilon = 0$ these relations reduce 
%with $f:=r-\ell$ and $c:=r\kappa - \ell \lambda$ to
%\begin{eqnarray}
%\sum_{k=2}^{L-2} h_{k} | \mu \rangle
%& = &  c \hat{n}_{1}\hat{n}_{2} \, | \mu \rangle  - [(f+c) \hat{n}_{2}  + c (\hat{n}_{1} -  \hat{n}_{2}) \hat{n}_{3} )] \, | \mu \rangle 
%\nonumber \\
%& &  - c \hat{n}_{L-1} \hat{n}_{L} \,| \mu \rangle
%+ [(f+c) \hat{n}_{L-1} + c \hat{n}_{L-2} (\hat{n}_{L} - \hat{n}_{L-1}) ] \,| \mu \rangle 
%\label{bulksum0} 
%\end{eqnarray}
%\begin{eqnarray}
%h^{-} | \mu \rangle
%& = &  - d^{-}_{1} (\hat{n}_{1} -  \hat{n}_{2}) \,  |\mu\rangle  
%+ (d^{-}_{1} - d^{-}_{2}) (\hat{n}_{1} -  \hat{n}_{2}) \hat{n}_{3} \,  |\mu\rangle
%\label{lb0} \\
%h^{+} | \mu \rangle
%& = & d^{+}_{1} (\hat{n}_{L} -  \hat{n}_{L-1}) \, | \mu \rangle
%- (d^{+}_{1} - d^{+}_{2}) \hat{n}_{L-2} (\hat{n}_{L} - \hat{n}_{L-1}) 
%\, | \mu \rangle
%\label{rb0} 
%\end{eqnarray}
%\begin{eqnarray}
%b^{-} | \mu \rangle  
%& = & [ - a^{-}_1 z + a^{-}_1 (1+z)\hat{n}_{1} + (a^{-}_1-a_2^{-} ) z \hat{n}_{2}  - (a^{-}_1-a_2^{-} ) (1+z) \hat{n}_{1} \hat{n}_{2}] \,  | \mu \rangle 
%\label{lres0} \\
%b^{+} | \mu \rangle  & = & \left[ a_1^{+} z -  (a_1^{+} - a_2^{+}) z\hat{n}_{L-1} 
%- a_1^{+} (1+z) \hat{n}_{L} 
%+ [a_1^{+} - a_2^{+}] (1+z) \hat{n}_{L-1} \hat{n}_{L} \right] | \mu \rangle .
%\label{rres0}
%\end{eqnarray}
\end{lemma}

\proof 
Using the relations \eqref{skmmu} and the definitions 
\eqref{fdef} and \eqref{cdef} one finds from \eqref{hkrev} for bulk terms with $2 \leq k \leq L-2$ after some lengthy but straightforward algebra 
\begin{eqnarray*}
\hat{\mu}^{-1} h_{k} \, | \mu \rangle
& = & \tilde{h}_k \, |s\rangle \\
& = &  r  (1+\kappa) (\hat{n}_{k+1} - \hat{n}_{k}) | s \rangle  \\
& &  + r (\epsilon - \kappa) (\hat{n}_{k+1} \hat{n}_{k+2} - \hat{n}_{k-1} \hat{n}_{k}) | s \rangle \\
& & + r (\epsilon + \kappa) (\hat{n}_{k} \hat{n}_{k+2} - \hat{n}_{k-1} \hat{n}_{k+1}) | s \rangle \\ 
& &  - 2 r \epsilon (\hat{n}_{k} \hat{n}_{k+1} \hat{n}_{k+2} - \hat{n}_{k-1} \hat{n}_{k} \hat{n}_{k+1}) | s \rangle  \\
& & - \ell  (1+\lambda) (\hat{n}_{k+1} - \hat{n}_{k}) | s \rangle  \\
& & 
- \ell  (\epsilon + \lambda) (\hat{n}_{k}  \hat{n}_{k+2} - \hat{n}_{k-1}  \hat{n}_{k+1}) | s \rangle \\
& & - \ell (\epsilon - \lambda) (\hat{n}_{k+1} \hat{n}_{k+2}
-  \hat{n}_{k-1} \hat{n}_{k} ) | s \rangle \\
& & 
+ 2 \ell \epsilon (\hat{n}_{k}  \hat{n}_{k+1} \hat{n}_{k+2} - \hat{n}_{k-1} \hat{n}_{k} \hat{n}_{k+1}) | s \rangle \\
% % % % % % % % % %
& = &  (f + c) (\hat{n}_{k+1} - \hat{n}_{k}) | s \rangle  \\
& & + (f \epsilon + c) (\hat{n}_{k} \hat{n}_{k+2} - \hat{n}_{k-1} \hat{n}_{k+1}) | s \rangle \\ 
& &  + (f \epsilon - c) (\hat{n}_{k+1} \hat{n}_{k+2} - \hat{n}_{k-1} \hat{n}_{k}) | s \rangle \\
& &  - ((f \epsilon + c)+(f \epsilon - c)) (\hat{n}_{k} \hat{n}_{k+1} \hat{n}_{k+2} - \hat{n}_{k-1} \hat{n}_{k} \hat{n}_{k+1}) | s \rangle
\end{eqnarray*}

%\red{
%\begin{eqnarray*}
%\hat{\mu}^{-1} h_{k} | \mu \rangle
%& = &  [(f+c) \hat{n}_{k+1} + c \hat{n}_{k} \hat{n}_{k+2} - c \hat{n}_{k+1} (\hat{n}_{k+2} + \hat{n}_{k}) ] \,| s \rangle  \\
%& &  - [(f+c) \hat{n}_{k})  + c \hat{n}_{k-1} \hat{n}_{k+1}
%- c \hat{n}_{k} (\hat{n}_{k+1} + \hat{n}_{k-1})] \, | s \rangle 
%\end{eqnarray*}
%}
Multiplying both sides by $\hat{\mu}$ and
rearranging terms turns this into
\begin{eqnarray*}
h_{k} | \mu \rangle & = &  (f + c) (\hat{n}_{k+1} - \hat{n}_{k}) | \mu \rangle  \\
& & + (f \epsilon + c) (
\hat{n}_{k} \hat{v}_{k+1} \hat{n}_{k+2} - \hat{n}_{k-1} \hat{v}_{k} \hat{n}_{k+1}) | \mu \rangle \\ 
& &  + (f \epsilon - c) (\hat{v}_{k} \hat{n}_{k+1} \hat{n}_{k+2}
- \hat{n}_{k-1} \hat{n}_{k} \hat{v}_{k+1}) | \mu \rangle \\
& = &  (f + c) (\hat{n}_{k+1} - \hat{n}_{k}) | \mu \rangle  \\
& & + (f \epsilon + c) (
\hat{n}_{k} \hat{v}_{k+1} \hat{n}_{k+2} - \hat{n}_{k-1} \hat{v}_{k} \hat{n}_{k+1}) | \mu \rangle \\ 
& &  + (f \epsilon - c) [(1-\hat{v}_{k} \hat{v}_{k+2})\hat{n}_{k+1} 
- (1-\hat{v}_{k-1} \hat{v}_{k+1})\hat{n}_{k} ] | \mu \rangle .
\end{eqnarray*}
Using the telescopic property of the sum one arrives at \eqref{bulksum}.

In the same fashion one finds from \eqref{musma}, \eqref{musmb},
 \eqref{muspa}, \eqref{muspb}
%%%% b-->b-a, d --> d-c, then a -->r^{-}_{1}/(1-\epsilon), b --> r^{-}_{2}, c --> \ell^{-}_{1} q/(1+\epsilon) d --> \ell^{-}_{2} q %%%%%%%%%%%%
\begin{eqnarray*}
\hat{\mu}^{-1} h^{-} | \mu \rangle 
& = & (r^{-}_{1} -\ell^{-}_1)(q_{-} \hat{n}_{2} - \hat{n}_{1}) \, |s\rangle  \\
& & -  (q_{-}-1) (r^{-}_{1} 
- \ell^{-}_{1} ) \hat{n}_{1}\hat{n}_{2} \, |s\rangle  \\
& & - q_{-}  \left[ r^{-}_{1}  - r^{-}_{2}  (1+\epsilon) 
- \ell^{-}_{1}  + \ell^{-}_{2} (1+\epsilon)  \right] \hat{n}_{2}\hat{n}_{3} \, |s\rangle\\
& & + \left[ r^{-}_{1} - r^{-}_{2} (1-\epsilon) 
- \ell^{-}_{1}  + \ell^{-}_{2} (1-\epsilon) \right] \hat{n}_{1} \hat{n}_{3} \, |s\rangle\\
& &  + \left[ r^{-}_{1} (q_{-}-1)  + r^{-}_{2} (1-\epsilon) - r^{-}_{2} q_{-} (1+\epsilon)  \right] \hat{n}_{1}\hat{n}_{2}\hat{n}_{3} \, |s\rangle\\
& &
- \left[\ell^{-}_{1}(q_{-}-1)
  + \ell^{-}_{2} (1-\epsilon) - \ell^{-}_{2} q_{-} (1+\epsilon) \right] \hat{n}_{1}\hat{n}_{2}\hat{n}_{3} \, |s\rangle
\end{eqnarray*}
% % % % % % % % % % % % % %
\begin{eqnarray*}
\hat{\mu}^{-1} h^{+} | \mu \rangle 
& = & - (r^{+}_{1} -\ell^{+}_1) (q_{+} \hat{n}_{L-1} - \hat{n}_{L}) \, |s\rangle \\
& & +(q_{+} - 1) (r^{+}_{1} -\ell^{+}_1) \hat{n}_{L-1}\hat{n}_{L} \, |s\rangle \\
& & + q_{+} [r^{+}_{1} - r^{+}_{2} (1+\epsilon) 
- \ell^{+}_1 + \ell^{+}_2 (1+\epsilon) ] \hat{n}_{L-2} \hat{n}_{L-1} \, |s\rangle\\
& & - [r^{+}_{1} - r^{+}_{2} (1-\epsilon)
- \ell^{+}_1 + \ell^{+}_2 (1-\epsilon) ]\hat{n}_{L-2}\hat{n}_{L} \, |s\rangle\\
& & - [r^{+}_{1} (q_{+} - 1) 
+ r^{+}_{2} (1-\epsilon)
- r^{+}_{2} q_{+} (1+\epsilon) ]\hat{n}_{L-2} \hat{n}_{L-1}\hat{n}_{L} \, |s\rangle\\
& & 
+ [\ell^{+}_1 (q_{+} - 1) + \ell^{+}_2 (1-\epsilon) 
- \ell^{+}_2 q_{+} (1+\epsilon) ] \hat{n}_{L-2}\hat{n}_{L-1}\hat{n}_{L} \, |s\rangle .
\end{eqnarray*}
By multiplying both sides with $\hat{\mu}$ and using the definitions 
\eqref{dpmidef} the
assertions \eqref{lb} and \eqref{rb} follow.

Finally, using \eqref{mus1} and \eqref{musL} one obtains
\begin{eqnarray*}
\hat{\mu}^{-1} b^{-}_{1} | \mu \rangle 
& = &  - (\alpha_1 z_{-} - \gamma_1 z_{-}) \, | s \rangle \\
& & + (\alpha_1 z_{-} - \gamma_1 z_{-}) (1+z_{-}^{-1}) \hat{n}_{1}  \, | s \rangle \\
& & +[\alpha_1 z_{-} - \gamma_1 z_{-} - (\alpha_2 z_{-} -  \gamma_2 z_{-} y) ] \hat{n}_{2} \\
& & - \left[(\alpha_1 z_{-} - \gamma_1 z_{-}) (1+z_{-}^{-1})- (\alpha_2 z_{-} -  \gamma_2 z_{-} y)  (1+y^{-1} z_{-}^{-1})\right] \hat{n}_{1}  
 \hat{n}_{2}\, | s \rangle \\
% % % % % % % % % % % % % % % % % % % % % % % % % % % % %
% % % % % % % % % % % % % % % % % % % % % % % % % % % % %
% % % % % % % % % % % % % % % % % % % % % % % % % % % % %
%& = & - (\alpha_1 z_{-} - \gamma_1 x) | s \rangle 
%\\ & & 
%+ (\alpha_1 z_{-}  - \gamma_1 x) (1+x^{-1}) \hat{n}_1 | s \rangle 
%\\ & & + [(\alpha_1 z_{-} - \gamma_1 x) - (\alpha_2 z_{-}   - \gamma_2 x y)]  \hat{n}_{2} | s \rangle 
%\\ & & 
%- [(\alpha_1 z_{-}  - \gamma_1 x)  (1+x^{-1}) - (\alpha_2 z_{-}   - \gamma_2 x y) 
%(1 + y^{-1} x^{-1})] \hat{n}_1 \hat{n}_{2} | s \rangle \\
% % % % % % % % % % % % % % % % % % % % % % % % % % % % %
% % % % % % % % % % % % % % % % % % % % % % % % % % % % %
% % % % % % % % % % % % % % % % % % % % % % % % % % % % %
\hat{\mu}^{-1} b^{+}_{L} | \mu \rangle 
%& = & (z^{-1} \hat{v}_{L-1} \hat{n}_{L}^{-} + y^{-1} z^{-1} \hat{n}_{L-1} \hat{n}_{L}^{-} - \hat{v}_{L}) (\delta_1 z_{+} \hat{v}_{L-1} + \delta_2 z_{+} \hat{n}_{L-1}) | s \rangle \\
%& & + (z \hat{v}_{L-1} \hat{v}_{L} + y z \hat{n}_{L-1} \hat{v}_{L} - \hat{n}_{L}) (\beta_1(\hat{v}_{L-1} + \beta_2 \hat{n}_{L-1}) | s \rangle \\
%& = & \delta_1 z_{+} (z^{-1} \hat{v}_{L-1} \hat{n}_{L}^{-} + y^{-1} z^{-1} \hat{n}_{L-1} \hat{n}_{L}^{-} - \hat{v}_{L}) \hat{v}_{L-1} | s \rangle \\
%& & + \delta_2 z_{+} (z^{-1} \hat{v}_{L-1} \hat{n}_{L}^{-} + y^{-1} z^{-1} \hat{n}_{L-1} \hat{n}_{L}^{-} - \hat{v}_{L}) \hat{n}_{L-1} | s \rangle \\
%& & + \beta_1 (z \hat{v}_{L-1} \hat{v}_{L} + y z \hat{n}_{L-1} \hat{v}_{L} - \hat{n}_{L})   \hat{v}_{L-1} | s \rangle \\
%& & + \beta_2 (z \hat{v}_{L-1} \hat{v}_{L} + y z \hat{n}_{L-1} \hat{v}_{L} - \hat{n}_{L}) \hat{n}_{L-1} | s \rangle \\
%& = & \delta_1 z_{+} (z^{-1} \hat{v}_{L-1} \hat{n}_{L}^{-} - \hat{v}_{L-1} \hat{v}_{L}) | s \rangle \\
%& & + \delta_2 z_{+} (y^{-1} z^{-1} \hat{n}_{L-1} \hat{n}_{L}^{-} 
%- \hat{n}_{L-1} \hat{v}_{L}) | s \rangle \\
%& & + \beta_1 (z \hat{v}_{L-1} \hat{v}_{L} - \hat{v}_{L-1} \hat{n}_{L}) | s \rangle \\
%& & + \beta_2 (y z \hat{n}_{L-1} \hat{v}_{L} - \hat{n}_{L-1} \hat{n}_{L}) | s \rangle \\
%& = & \left[ (\beta_1 z - \delta_1 z_{+}) \hat{v}_{L-1} \hat{v}_{L}
%+ (\delta_1 z_{+} z^{-1} - \beta_1) \hat{v}_{L-1} \hat{n}_{L} \right] | s \rangle 
%\\ & &
%+ \left[ (\beta_2 y z - \delta_2 z_{+}) \hat{n}_{L-1} \hat{v}_{L} 
%+ (\delta_2 z_{+} y^{-1} z^{-1} - \beta_2) \hat{n}_{L-1} \hat{n}_{L} \right] | s \rangle \\
& = & (\beta_1 z_{+} - \delta_1 z_{+}) | s \rangle  
\\ & & 
-  [(\beta_1 z_{+} - \delta_1 z_{+}) - (\beta_2 yz_{+} - \delta_2 z_{+})] \hat{n}_{L-1}  | s \rangle
\\ & & 
- (\beta_1 z_{+} - \delta_1 z_{+}) (1+z_{+}^{-1}) \hat{n}_{L}  | s \rangle
\\ & &
+ [(\beta_1 z_{+} - \delta_1 z_{+}) (1+z_{+}^{-1})  - (\beta_2 y z_{+} - \delta_2 z_{+})(1 + y^{-1} z_{+}^{-1} )] \hat{n}_{L-1} \hat{n}_{L} | s \rangle  .
\end{eqnarray*}
By multiplying both sides with $\hat{\mu}$ and using the definitions 
\eqref{c1mdef} - \eqref{c2pdef} one arrives at \eqref{lres} and \eqref{rres}.
\qed

\paragraph{(iii) Proof of invariance of the Ising measure:}
Recall that $z_{\pm} = q_\pm z$.
From Lemma \ref{Lmm:Hmu} one obtains
\begin{eqnarray}
H \, |\mu\rangle
& = & \left[a_1^{-} (1+z_{-})  - d^{-}_{1}\right] \hat{n}_{1}  \,  |\mu\rangle 
\label{cH1m} \\
& &  + \left[(a_1^{-} - a_2^{-}) z_{-} 
+ q_{-} d^{-}_{1} - (f + c)\right] \hat{n}_{2} \,  |\mu\rangle 
\label{cH2m}\\
& & - \left[
a_1^{-}  (1+z_{-}) - a_2^{-} 
(z_{-} + y^{-1}) + (q_{-}-1) d^{-}_{1}
+ (f \epsilon - c)\right]  \hat{n}_1 \hat{n}_{2} \,| \mu \rangle \label{cH12m} \\
& & - \left[- d^{-}_{1} +d^{-}_{2}(1-\epsilon) + (f \epsilon + c)\right] \hat{n}_{1} \hat{n}_{3} \,  |\mu\rangle 
\label{cH13m}\\
& &  - \left[d^{-}_{1} q_{-} - d^{-}_{2}q_{-} (1+\epsilon) + (f \epsilon - c)\right]  \hat{n}_{2}\hat{n}_{3} \,  |\mu\rangle 
\label{cH23m} \\
& & + \left[2 f \epsilon  + d^{-}_{1}(q_{-}-1)  + d^{-}_{2} (1-\epsilon)
- d^{-}_{2} q_{-} (1+\epsilon) \right] \hat{n}_{1}\hat{n}_{2}\hat{n}_{3} \,  |\mu\rangle 
\label{cH123m}\\[2mm]
% % % % % % % % % % % % % % % %
% % % % % % % % % % % % % % % %
% % % % % % % % % % % % % % % %
% % % % % % % % % % % % % % % %
& & - \left[a_1^{+} (1+z_{+}) - d^{+}_{1}\right] \hat{n}_{L} \,  |\mu\rangle 
\label{cH1p}\\
& & - \left[ (a_1^{+} - a_2^{+}) z_{+}\hat{n}_{L-1} 
+ q_{+} d^{+}_{1} \hat{n}_{L-1}
- (f + c) \right] \hat{n}_{L-1}
\label{cH2p}\\
& & + \left[  
a_1^{+} (1+z_{+})  - a_2^{+} (z_{+} + y^{-1}) + (q_{+} - 1) d^{+}_{1} + (f \epsilon - c) \right] 
\hat{n}_{L-1} \hat{n}_{L}\, | \mu \rangle  
\label{cH12p}\\
& &  + \left[- d^{+}_{1}
+ d^{+}_{2} (1-\epsilon) + (f \epsilon + c) \right] \hat{n}_{L-2}\hat{n}_{L} \, | \mu \rangle 
\label{cH13p} \\
& & + \left[d^{+}_{1} q_{+} 
- d^{+}_{2} q_{+} (1+\epsilon) + (f \epsilon - c) \right]  \hat{n}_{L-2}\hat{n}_{L-1} 
\, | \mu \rangle 
\label{cH23p}\\
& & - \left[ 
 2 f \epsilon  + d^{+}_{1} (q_{+} - 1) + d^{+}_{2} (1-\epsilon) - d^{+}_{2} q_{+} (1+\epsilon) \right] \hat{n}_{L-2} \hat{n}_{L-1} \hat{n}_{L}) \, |\mu \rangle 
\label{cH123p} \\[2mm]
% % % % % % % % % % % % % % % %
% % % % % % % % % % % % % % % %
% % % % % % % % % % % % % % % %
% % % % % % % % % % % % % % % %
& & + (a_1^{+} z_{+}  - a_1^{-} z_{-}) \,  |\mu\rangle \label{cHpm}
\end{eqnarray}
The right hand side of the equation is generically nonzero but allows for proving the theorem by matching 
parameters for cancellation
of the ``unwanted'' terms inside the braces.
Notice that invariance holds if and only if all terms
\eqref{cH1m} - \eqref{cHpm} vanish. The proof
is facilitated by the fact the cancellation conditions
for the left boundary terms \eqref{cH1m} - \eqref{cH123m} are, up to the sign, identical to \eqref{cH1p} - \eqref{cH123p} arising from the
right boundary.
We first prove the cancellation for $\epsilon \neq 0$.\\
(i) The conditions \eqref{apm1} and \eqref{cH1}
are immediate consequences of \eqref{cH1m}, \eqref{cH1p}, and \eqref{cHpm}.\\
(ii) For $\epsilon \neq 0$ the choice of rates \eqref{dpm1} and \eqref{dpm2} is equivalent to
\begin{eqnarray}
(f \epsilon + c) & = & d^{-}_{1}  - d^{-}_{2} (1-\epsilon) 
\,=\, d^{+}_{1}  - d^{+}_{2}(1-\epsilon)
\label{c1} \\
(f \epsilon - c) & = & - q_{-} [d^{-}_{1}- d^{-}_{2} (1+\epsilon)]
\,=\, - q_{+} [d^{+}_{1} - d^{+}_{2}(1+\epsilon)]
\label{c2}
\end{eqnarray}
which proves cancellation of the terms \eqref{cH13m} - \eqref{cH123m}
and \eqref{cH13p} - \eqref{cH123p}.\\
(iii) To complete the cancellations it remains to prove that
\begin{eqnarray}
0 & = &  (a_1^{\pm} - a_2^{\pm}) z_{\pm} 
+ q_{\pm} d^{\pm}_{1} - (f + c) 
\label{cH2pm}\\
0 & = &  
a_1^{\pm}  (1+z_{\pm}) - a_2^{\pm} 
(z_{\pm} + y^{-1}) + (q_{\pm}-1) d^{\pm}_{1}
+ (f \epsilon - c)  \label{cH12pm}
\end{eqnarray}
Observing that
\begin{equation}
(f + c) = \frac{(f \epsilon + c) (1+\epsilon) + (f \epsilon - c) (1-\epsilon)}{2 \epsilon}.
\end{equation}
one finds
\begin{eqnarray*}
d^{\pm}_{1} & = & 
(f + c) + (q_{\pm}^{-1} - 1) \frac{(f \epsilon - c) (1-\epsilon)}{2 \epsilon}\\
q_{\pm} d^{\pm}_{1} & = & 
(f + c) + (q_{\pm} - 1) \frac{(f \epsilon + c) (1+\epsilon)}{2 \epsilon}
\end{eqnarray*}
With \eqref{dpm1} - \eqref{apm2} it follows that
\eqref{cH2pm} and \eqref{cH12pm} are
equivalent to cancellation of the terms \eqref{cH2m}, \eqref{cH12m},
\eqref{cH2p}, and \eqref{cH12p}.
Hence $H | \mu \rangle = 0$ if and only if the
conditions \eqref{dpm1} - \eqref{cH1} are satisfied.\\[4mm]

(2) The special case $\epsilon=0$ with $c\neq 0$ is obtained
as the limit $\epsilon\to 0$ of the generic case 
$0<|\epsilon|<1$ by observing that
solving the conditions \eqref{cH4} and \eqref{cH7}
for $z_{\pm}$ and expanding in powers of $\epsilon$ yields
$q_{\pm} = 1 - 2\epsilon z/(1+z) + O(\epsilon^2)$.\\[4mm]

The proof for $\epsilon=c=0$ is slightly different. As above,
the contributions \eqref{cH123m} and \eqref{cH123p} can be ignored
since they automatically cancel if \eqref{cH13m}, \eqref{cH23m}, 
\eqref{cH13p} and \eqref{cH23p} cancel. The contribution \eqref{cHpm}
yields the compatibility condition \eqref{cH1}.  The remaining
five pairs of equalities \eqref{cH1m} - \eqref{cH23m} 
and \eqref{cH1p} - \eqref{cH23p}  reduce to the equalities
\begin{eqnarray}
0 & = & a_1^{\pm} (1+z_{\pm})  - d^{\pm}_{1}
\label{cH1pm00} \\
0 & = &  (a_1^{\pm} - a_2^{\pm}) z_{\pm} 
+ q_{\pm} d^{\pm}_{1} - f
\label{cH2pm00}\\
0 & = & 
a_1^{\pm}  (1+z_{\pm}) - a_2^{\pm} 
(1+ z_{\pm}) + (q_{\pm}-1) d^{\pm}_{1}  \label{cH12pm00} \\
0 & = &  - d^{\pm}_{1} +d^{\pm}_{2} 
\label{cH13pm00}\\
0 & = &  d^{\pm}_{1} q_{\pm} - d^{\pm}_{2}q_{\pm}  
\label{cH23pm00} 
\end{eqnarray}
which are necessary and sufficient for invariance.

The last two equations \eqref{cH13pm00} and \eqref{cH23pm00}
are simultaneously equivalent to $d^{\pm}_{2} = d^{\pm}_{1}$.
The first equation \eqref{cH1pm00} yields $a_1^{\pm} (1+z_{\pm})  = d^{\pm}_{1}$. Inserting this into \eqref{cH2pm00} and \eqref{cH12pm00}
and solving for $a_2^{\pm}$ 
yields
\begin{eqnarray}
 d^{\pm}_{1} & = &  f 
\frac{1+ z_{\pm}^{-1}}{1 + z^{-1}}  \\
a_2^{\pm} (1+ z_{\pm}) & = &  q_{\pm} d^{\pm}_{1} 
\end{eqnarray}
These equations are equivalent to the conditions \eqref{d1pm00},
 \eqref{d2pm00},  \eqref{a1pm00}, and \eqref{a2pm00}
for arbitrary $z_{\pm}\in\R^+$. \qed

\subsubsection{Proof of Theorem \ref{Theo:gKLSopenrev}}

With Proposition \ref{Prop:gst2} the proof is entirely analogous to the proof of Theorem \ref{Theo:gKLSrev}
and omitted here.

\section*{Acknowledgements}
This work is financially supported 
by FCT (Portugal) through project UIDB/04459/2020, 
doi 10-54499/UIDP/04459/2020,
and by the FCT Grants  2020.03953.CEECIND and 2022.09232.PTDC.  
N. Ngoc gratefully acknowledges the financial support of CAPES and CNPq during his 
studies at the Doctorate Program in Statistics at the University of S\~ao Paulo.

\appendix

\section{Expectation values for Ising measures}
\label{App:Ising}

We summarize the gist of the transfer matrix approach to the computation of 
expectation values for the Ising measure, see e.g. \cite{Baxt82} for a more detailed exposition. 

\paragraph{Linear algebra conventions:} 
We denote by $\bra{0} \equiv (1,0)$, $\bra{1} \equiv (0,1)$ the canonical basis vectors of $\R^2$ expressed as row vectors.
Introducing for any two-dimensional
row vector $\bra{a}$ with components $a_0,a_1\in\R$ the 
column vector $\ket{a} =\bra{a}^T$, the scalar product $<a,b>=a_0 b_0 + a_1 b_1$ 
of two vectors in $\R^2$ can be expressed as a matrix product, denoted by 
$\inprod{a}{b}$ of the row vector $\bra{a}$ (understood as a rectangular 
$1\times 2$ matrix) with the column vector $\ket{b}$ (understood as a rectangular 
$2\times 1$ matrix). The Kronecker product $\ket{a}\otimes\bra{b}$ is the $2\times 2$ 
matrix $C$ with elements $C_{ij} = a_{i} b_{j}$ and is denoted by
$C=\ket{a}\bra{b}$ without the Kronecker product symbol $\otimes$.
We also introduce the projectors
\begin{equation}
\mathsf{v} = \left(\ba{cc} 1 & 0 \\ 0 & 0 \ea \right) = \ket{0}\bra{0}, \quad 
\mathsf{n} = \left(\ba{cc} 0 & 0 \\ 0 & 1 \ea \right) = \ket{1}\bra{1} 
\end{equation}
expressed as Kronecker products of the canonical basis vectors $\bra{\eta}$ 
of $\R^2$.  The two-dimensional unit matrix is denoted by 
$\mathds{1}$ and can be expressed as the sums 
\begin{eqnarray}
\mathds{1} & = & \ket{0}\bra{0} + \ket{1}\bra{1} 
\label{unit1} 
\end{eqnarray}
of Kronecker products.

\paragraph{Transfer matrix:} In terms of the local Boltzmann weights
\begin{equation}
\mathsf{T}_{\eta,\eta'} :=\mathrm{e}^{-\beta [J \eta \eta' - \phi (\eta+\eta')/2]}.
\end{equation}
the full Boltzmann weight \eqref{IsingBW} of the periodic
generalized KLS model can be written 
\begin{equation}
\pi_{L}(\bfeta) = \mathsf{T}_{\eta_1,\eta_2} \mathsf{T}_{\eta_2,\eta_3} \dots
\mathsf{T}_{\eta_{L-1},\eta_L} \mathsf{T}_{\eta_L,\eta_1}.
\end{equation}
The idea of the transfer matrix approach is to consider the local Boltzmann weights as matrix elements of a matrix
\begin{equation}
\mathsf{T} := \left(\ba{cc} 1 & \sqrt{z} \\ \sqrt{z} & y z \ea \right),
\label{transfermatrix}
\end{equation} 
%\red{
%\begin{eqnarray*}
%\mathsf{T} |v> 
%& = &
%\left(\ba{cc} 1 & \sqrt{z} \\ \sqrt{z} & y z \ea \right)
%\left(\ba{c} \lambda - yz  \\ \sqrt{z} \ea \right) \\
%& = & \left(\ba{c} \lambda - yz + z \\ \sqrt{z} \lambda \ea \right) \\
%%& = & \left(\ba{c} \lambda^2 - yz \lambda 
%%- [\lambda^2 - (1+yz) \lambda + yz - z] \\ \sqrt{z} \lambda \ea \right) \\
%& = & \lambda \left(\ba{c} \lambda - yz  \\ \sqrt{z} \ea \right) 
%\end{eqnarray*} 
%\begin{eqnarray*}
%\mathsf{T}^2 & = & \left(\ba{cc} 1 & \sqrt{z} \\ \sqrt{z} & y z \ea \right)\left(\ba{cc} 1 & \sqrt{z} \\ \sqrt{z} & y z \ea \right) \\
%& = & \left(\ba{cc} 1 + z & \sqrt{z} (1+yz) \\ \sqrt{z} (1+yz) & z+ y^2 z^2 \ea \right)
%\end{eqnarray*}
%}
called transfer matrix.
The eigenvalues 
%\red{
%\begin{equation}
%0 = yz - z - (1+yz)\lambda + \lambda^2
%\end{equation}
%\begin{equation}
%\lambda = yz - z - yz \lambda + \lambda^2
%\end{equation}
%\begin{equation}
%\lambda^2 - 2yz \lambda + y^2z^2 = z(1-y+y^2z) + (1-yz)\lambda
%\end{equation}
%\begin{eqnarray*}
%a + z & = & a \lambda \\
%a + y\sqrt{z}  & = & \lambda \\
%\lambda - y\sqrt{z} + z & = & \lambda^2 - y\sqrt{z} \lambda  \\
%\end{eqnarray*}
%}
are
\begin{equation}
\lambda_{j} = \frac{1}{2} \left(1+y z + (-1)^{j} \sqrt{(1-y z)^2+4 z}\right), \quad j\in\{0,1\} .
\label{eigenvalues}
\end{equation}
%\red{
%\begin{eqnarray*}
%(\lambda_j-yz)^2 + z
%& = & z+ \frac{1}{2} \left[(1-yz)^2 + 2z + (-1)^j (1-yz) \sqrt{1-yz)^2 + 4z}\right] \\
%& = & z+\frac{1}{2} \left[2z + (1-yz)\left(1-yz  + (-1)^j \sqrt{(1-yz)^2 + 4z}\right) \right] \\
%& = & 2z + (1-yz) (\lambda_j-yz) \\
%%%%%%%%%%%%%%%%%%%%%%
%((\lambda_0-yz)^2 + z) \lambda_1 
%& = & z [2 \lambda_1 + (1-yz) (y-1 -y \lambda_1 )] \\
%& = & z [(2-y+y^2z) \lambda_1  + (1-yz) (y-1)] \\
%((\lambda_1-yz)^2 + z) \lambda_0
%& = & z [(2-y+y^2z) \lambda_0  + (1-yz) (y-1)] \\[4mm]
%\Sigma & = & z [(2-y+y^2z)(1+yz)  + 2 (1-yz) (y-1)] \\
%\Sigma & = & z y [(1-yz)^2+4z] \\[4mm]
%((\lambda_0-yz)^2 + z)((\lambda_1-yz)^2 + z) 
%& = & [2z + (1-yz) (\lambda_0-yz)][2z + (1-yz) (\lambda_1-yz)] \\
%& = & 4z^2 + 2z (1-yz) (1+yz-2yz) + (1-yz)^2  (\lambda_1-yz)(\lambda_0-yz) \\
%& = & z [4z  + (1-yz)^2]   \\
%\end{eqnarray*}
%}
and the normalized eigenvectors 
$\bra{v^{(j)}} := (v^{(j)}_0, v^{(j)}_1)$
with components
\begin{equation}
v^{(j)}_0 = \frac{\lambda_{j} - y z}{\sqrt{(\lambda_{j} - y z)^2+z}}, \quad
v^{(j)}_1 = \frac{\sqrt{z}}{\sqrt{(\lambda_{j} - y z)^2+z}}.
\label{eigenvectors}
\end{equation}
%\red{
%\begin{eqnarray*}
%Y_{1} & = & (v^{(0)}_1)^2 \lambda_{0} + (v^{(1)}_1)^2 \lambda_{1} \\
%& = &  \frac{z}{(\lambda_{0} - y z)^2+z} \lambda_{0} 
%+ \frac{z}{(\lambda_{1} - y z)^2+z}  \lambda_{1} \\
%& = &  z \left( \frac{[(\lambda_{1} - y z)^2+z]\lambda_{0} 
%+[(\lambda_{0} - y z)^2+z]\lambda_{1}}{[(\lambda_{0} - y z)^2+z][(\lambda_{1} - y z)^2+z]} 
% \right) \\
%\end{eqnarray*}
%}
form an orthogonal basis of $\R^2$, i.e.,
$\inprod{v^{(j)}}{v^{(j)}}=1$ and $\inprod{v^{(j)}}{v^{(1-j)}}=0$.
Orthogonality also implies the eigenvector decomposition 
\begin{equation}
\mathsf{T}^n = \lambda_{0}^n \ket{v^{(0)}}\bra{v^{(0)}} + \lambda_{1}^n \ket{v^{(1)}}\bra{v^{(1)}}, \quad n\in\N_0
\label{Teigen}
\end{equation}
which for $n=0$ yields the alternative representation
\begin{equation}
\mathds{1} = \ket{v^{(0)}}\bra{v^{(0)}} +\ket{v^{(1)}}\bra{v^{(1)}}
\label{unit2}
\end{equation}
of the unit matrix.

\paragraph{Expectation values for periodic boundary conditions:} 
By definition, the partition function \eqref{Zper} is given by
\begin{equation}
Z_{L} = \sum_{\eta_1\in\S} \sum_{\eta_2\in\S} \dots 
\sum_{\eta_{L-1} \in\S} \sum_{\eta_L\in\S} \mathsf{T}_{\eta_1,\eta_2} \mathsf{T}_{\eta_2,\eta_3} \dots
\mathsf{T}_{\eta_{L-1},\eta_L} \mathsf{T}_{\eta_L,\eta_1}
\end{equation}
so that, using $\mathsf{T}_{\eta,\eta'}
= \bra{\eta} \mathsf{T} \ket{\eta'} $ and the representation \eqref{unit1} of the unit matrix
\begin{eqnarray}
Z_{L} & = & \sum_{\eta_1\in\S} \sum_{\eta_2\in\S} \dots 
\sum_{\eta_{L-1} \in\S} \sum_{\eta_L\in\S}  \bra{\eta_1} \mathsf{T} \ket{\eta_2} 
\bra{\eta_2} \mathsf{T} \ket{\eta_3} \dots \bra{\eta_{L-1}} \mathsf{T} \ket{\eta_L} 
\bra{\eta_{L}} \mathsf{T} \ket{\eta_1} \nonumber \\
& = & \sum_{\eta_1\in\S} \bra{\eta_1} \mathsf{T}^L \ket{\eta_1} \, = \, \mathrm{Tr} \left(\mathsf{T}^{L} \right) \, = \, \lambda_{0}^L + \lambda_{1}^L.
\end{eqnarray}
Observing that
\begin{equation}
\sum_{\eta\in\S} \eta \ket{\eta}\bra{\eta} = \ket{1}\bra{1} = \mathsf{n}
\end{equation}
it follows that joint expectations of the occupation variables $\eta_k$ 
are given by the normalized trace
\begin{equation}
\exval{\eta_{k_1}\eta_{k_2}\cdots\eta_{k_n}}_{\mu} = 
\frac{1}{Z_L} \mathrm{Tr} \left(\mathsf{T}^{k_1} \mathsf{n} \mathsf{T}^{k_2-k_1} \mathsf{n} \cdots \mathsf{T}^{k_n-k_{n-1}} 
\mathsf{n} \mathsf{T}^{L-k_n} \right).
\label{exvalIsing}
\end{equation}

To compute these expectations it is convenient to introduce for $n\in\N$ the quantities
\begin{equation}
Y_n := \bra{1} \mathsf{T}^n \ket{1} = (v^{(0)}_1)^2 \lambda_{0}^{n} + (v^{(1)}_1)^2 \lambda_{1}^{n} 
\label{Ydef}
\end{equation}
where the equality follows from the eigenvalue decomposition \eqref{Teigen}. 
Translation invariance and the cyclic property of the trace
thus yields
\begin{eqnarray}
\rho(L) & = & \frac{1}{Z_L} \mathrm{Tr} \left(\mathsf{n} \mathsf{T}^{L} \right) =  \frac{Y_L}{Z_L}
\label{rhoIsing1} \\
\exval{\eta_k\eta_{k+r}}_L & = & \frac{1}{Z_L} \mathrm{Tr} \left(\mathsf{n} \mathsf{T}^{r} \mathsf{n} \mathsf{T}^{L-r}\right) \, = \,  \frac{Y_r Y_{L-r}}{Z_L} 
\label{GIsing1} \\
\exval{\eta_k\eta_{k+1}\eta_{k+2}}_L & = & \frac{1}{Z_L} \mathrm{Tr} \left(\mathsf{n} \mathsf{T} \mathsf{n} \mathsf{T} \mathsf{n} \mathsf{T}^{L-2}\right) \, = \,  \frac{Y_1^2 Y_{L-2}}{Z_L} 
\label{HIsing1} 
\end{eqnarray}
for all $k\in\T_L$.
It follows that
\begin{eqnarray}
\exval{\eta_k \bar{\eta}_{k+1}}_L & = & \frac{1}{Z_L} \mathrm{Tr} \left(\mathsf{n} \mathsf{T} \mathsf{v} \mathsf{T}^{L-1}\right) \, = \,  \frac{Y_L - Y_1 Y_{L-1}}{Z_L} 
\, = \, \frac{1}{Z_L} \mathrm{Tr} \left(\mathsf{v} \mathsf{T} \mathsf{n} \mathsf{T}^{L-1}\right) 
\label{GIsing2} \\
\exval{\eta_k \bar{\eta}_{k+1} \eta_{k+2}}_L & = & \frac{1}{Z_L} \mathrm{Tr} \left(\mathsf{n} \mathsf{T} \mathsf{v} \mathsf{T} \mathsf{n} \mathsf{T}^{L-2}\right)
\, = \, \frac{(Y_2 -
Y_1^2) Y_{L-2}}{Z_L} .
\label{HIsing21} \\
\exval{\eta_k \eta_{k+1} \bar{\eta}_{k+2}}_L  & = & \frac{1}{Z_L} \mathrm{Tr} \left(\mathsf{n} \mathsf{T} \mathsf{n} \mathsf{T} \mathsf{v} \mathsf{T}^{L-2}\right) 
\, = \, \frac{Y_1 (Y_{L-1} -
Y_1 Y_{L-2})}{Z_L} = \frac{1}{Z_L} \mathrm{Tr} \left(\mathsf{v} \mathsf{T} \mathsf{n} \mathsf{T} \mathsf{n} \mathsf{T}^{L-2}\right).
\label{HIsing22}  
\end{eqnarray}

Since $Y_1 = yz $ and $Y_2 = z(1 + y^2 z)$ we arrive at the explicit expressions
\begin{eqnarray}
\exval{\eta_k \bar{\eta}_{k+1}}_L & = & \exval{\bar{\eta}_k \eta_{k+1}}_L
\, = \,  \left(1 - y z \frac{Y_{L-1}}{Y_L} \right) \rho(L) 
\label{GIsing3} \\
\exval{\eta_k \bar{\eta}_{k+1} \eta_{k+2}}_L & = &  z 
\frac{ Y_{L-2}}{Y_L} \rho(L)
\label{HIsing31} \\
\exval{\eta_k \eta_{k+1} \bar{\eta}_{k+2}}_L  & = & \exval{\bar{\eta}_k \eta_{k+1} \eta_{k+2}}_L
\, = \, y z \left(\frac{Y_{L-1}}{Y_L}  -
y z \frac{Y_{L-2}}{Y_L} \right) \rho(L)
\label{HIsing32}  
\end{eqnarray}
for the expectations that appear in the stationary current \eqref{jLdef}.

\paragraph{Thermodynamic limit:} 
To study the thermodynamic limit it is convenient to introduce the ratios
\begin{equation}
b := \frac{v^{(1)}_{1}}{v^{(0)}_{1}}, \quad c := \frac{\lambda_1}{\lambda_0}
\end{equation}
and take note of the fact that $c<1$ since by definition $\lambda_1<\lambda_0$.
In terms of these constants the density reads
\begin{equation}
\rho(L) = (v^{(0)}_1)^2 \frac{1 + b^2 c^{L} }{1 + c^{L} }
\end{equation}
so that
\begin{equation}
\rho = \lim_{L\to\infty} \rho(L) = (v^{(0)}_1)^2
\end{equation}
follows immediately. In a similar fashion one obtains
\begin{eqnarray}
\lim_{L\to\infty} \exval{\eta_k \bar{\eta}_{k+1}}_L & = & 
\lim_{L\to\infty} \exval{\bar{\eta}_k \eta_{k+1}}_L
\, = \,  \left(1 - \frac{y z}{\lambda_1} \right) \rho 
\\
\lim_{L\to\infty}  \exval{\eta_k \bar{\eta}_{k+1} \eta_{k+2}}_L & = &  \frac{z}{\lambda_1^2} \rho
 \\
\lim_{L\to\infty}  \exval{\eta_k \eta_{k+1} \bar{\eta}_{k+2}}_L  & = & \lim_{L\to\infty}   \exval{\bar{\eta}_k \eta_{k+1} \eta_{k+2}}_L
\, = \, \frac{y z}{\lambda_1} \left(1  -
\frac{y z}{\lambda_1} \right) \rho.
\end{eqnarray}
Finite-size corrections to these limiting values are exponentially small
in $L$ as the appearance of the term $c^L$ in the expressions for finite
size demonstrates.

\paragraph{Ising model with boundary fields:} 
The Boltzmann weight of the one-dimensional Ising model with boundary fields
is given by
\begin{equation}
\pi_{L}(\boldsymbol{\eta}) 
= \mathrm{e}^{-\beta \left[J \sum_{k=1}^{L-1}\eta_k\eta_{k+1} - \varphi \sum_{k=2}^{L-1}\eta_k - \varphi_{-} \eta_1  - \varphi_{+} \eta_L\right]}
%= \mathrm{e}^{-\beta \left[J \sum_{i=1}^{L-1}\eta_i\eta_{i+1} - \frac{\varphi}{2} \sum_{i=1}^{L-1}(\eta_i+\eta_{i+1}) - (\varphi_{-}-\frac{\varphi}{2}) \eta_1  - (\varphi_{+}-\frac{\varphi}{2}) \eta_L\right]}
\label{IsingBWfree}
\end{equation}
The computation of expectation values with the transfer matrix approach is analogous. However, since the interaction term
between sites $L$ and $1$ is replaced by boundary terms, the trace
operation in the computation of expectations is replaced by a scalar product
with the vectors 
\begin{equation}
\bra{s_{-}} := \bra{0} + \mathrm{e}^{\beta  (\varphi_{-}-\frac{\varphi}{2})} \bra{1}, \quad \ket{s_{+}} := \ket{0} + \mathrm{e}^{\beta  (\varphi_{+}-\frac{\varphi}{2})} \ket{1}.
\end{equation} 
Specifically, for the partition function this yields
\begin{eqnarray}
Z_L & = & \sum_{\eta_1\in\S} \sum_{\eta_2\in\S} \dots 
\sum_{\eta_{L-1} \in\S} \sum_{\eta_L\in\S}  
\mathrm{e}^{\beta \left[ (\varphi_{-}-\frac{\varphi}{2}) \eta_1  + (\varphi_{+}-\frac{\varphi}{2}) \eta_L\right]}
\bra{\eta_1} \mathsf{T} \ket{\eta_2} 
\bra{\eta_2} \mathsf{T} \ket{\eta_3} \dots \bra{\eta_{L-1}} \mathsf{T} \ket{\eta_L} 
 \nonumber \\
& = & \sum_{\eta_1\in\S} \sum_{\eta_L\in\S} 
\mathrm{e}^{\beta \left[ (\varphi_{-}-\frac{\varphi}{2}) \eta_1  + (\varphi_{+}-\frac{\varphi}{2}) \eta_L\right]}
\bra{\eta_1} \mathsf{T}^{L-1} \ket{\eta_L} 
\, = \, \bra{s_{-}} \mathsf{T}^{L-1} \ket{s_{+}}.
\end{eqnarray}
With the coefficients
\begin{equation}
a^{(j)}_{\pm} := v^{(j)}_0 + \mathrm{e}^{\beta  (\varphi_{\pm}-\frac{\varphi}{2})} v^{(j)}_1 , \quad
a_{\pm} := \frac{a^{(1)}_{\pm}}{a^{(0)}_{\pm}}
%= \frac{v^{(1)}_0 + \mathrm{e}^{\beta  (\varphi_{\pm}-\frac{\varphi}{2})} v^{(1)}_1 }{ v^{(0)}_0 + \mathrm{e}^{\beta  (\varphi_{\pm}-\frac{\varphi}{2})} v^{(0)}_1 }
\end{equation}
The eigenvalue decomposition \eqref{Teigen} thus yields
\begin{equation}
Z_L =  a_{+}^{(0)} a_{-}^{(0)} \lambda_{0}^{L-1}  + 
a_{+}^{(1)} a_{-}^{(1)} \lambda_{1}^{L-1} 
= a_{+}^{(0)} a_{-}^{(0)} \lambda_{0}^{L-1}
\left(1+a_{+}a_{-} c^{L-1}\right).
\end{equation}
and after some straightforward computation the stationary density profile \eqref{rhokldef}
\begin{eqnarray}
\rho_k(L) 
& = & \frac{1}{Z_L} \bra{s_{-}} \mathsf{T}^{k-1} \mathsf{n} \mathsf{T}^{L-k} \ket{s_{+}} \\
& = & \left\{
\begin{array}{ll}
\displaystyle  \frac{\mathrm{e}^{\beta  (\varphi_{-}-\frac{\varphi}{2})} v_{1}^{(0)} \left(1+a_{+}bc^{L-1}\right) }{ a_{-}^{(0)} 
\left(1+a_{+}a_{-} c^{L-1}\right)} & k=1 \\[6mm]
\displaystyle \rho \, \frac{(1 + a_{+} b  c^{L-k})(1 + a_{-} b c^{k-1}) 
}{ 1 + a_{+} a_{-} c^{L-1}} & 2 \leq k \leq L-1 \\[6mm]
\displaystyle  \frac{\mathrm{e}^{\beta  (\varphi_{+}-\frac{\varphi}{2})} v_{1}^{(0)} \left(1+a_{-}bc^{L-1}\right) }{ a_{+}^{(0)} 
\left(1+a_{+}a_{-} c^{L-1}\right)} & k=L.
\end{array}
\right.
\label{densityprofile}
\end{eqnarray}
For $k=[xL]$ with $0<x<1$ the bulk density profile 
$\rho_{[xL]}(L)$ in the thermodynamic limit $L\to\infty$ is equal to the
constant density $\rho$ of the periodic system in the thermodynamic
limit. With similar computations one sees that this is true
also for all bulk correlations of finite order and consequently for the current. Notice also that
in the thermodynamic limit there are generically
boundary discontinuities, i.e., $\lim_{L\to\infty}
\rho_{1}(L) \neq \rho$ and $\lim_{L\to\infty}
\rho_{L}(L) \neq \rho$, unless $\mathrm{e}^{\beta  (\varphi_{-}-\frac{\varphi}{2})} = \mathrm{e}^{\beta  (\varphi_{+}-\frac{\varphi}{2})} = v^{(0)}_1/v^{(0)}_0$, i.e., unless the boundary
chemical potentials are chosen such that $\bra{s_{\pm}}$ are eigenvectors of the transfer
matrix with eigenvalue $\lambda_0$. In the expression \eqref{densityprofile} this is seen by
noting that with this choice of boundary fugacities
one has by orthogonality and normalization of the
eigenvectors $\bra{s_{\pm}}$ that $a^{(0)}_{\pm} =1/v^{(0)}_0$ and $a^{(1)}_{\pm} =0$. This implies $a_{\pm} =0$ 
and therefore $\rho_{k}(L) = \rho$ for all $L\geq 2$.

\end{document}